
\def\LaTeX{%
  \let\Begin\begin
  \let\End\end
  \let\salta\relax
  \let\finqui\relax
  \let\futuro\relax}

\def\UK{\def\our{our}\let\sz s}
\def\USA{\def\our{or}\let\sz z}

\UK



\LaTeX

\USA


\salta

\documentclass[twoside,12pt]{article}
\setlength{\textheight}{24cm}
\setlength{\textwidth}{16cm}
\setlength{\oddsidemargin}{2mm}
\setlength{\evensidemargin}{2mm}
\setlength{\topmargin}{-15mm}
\parskip2mm


\usepackage[usenames,dvipsnames]{color}
\usepackage{amsmath}
\usepackage{amsthm}
\usepackage{amssymb,bbm}
\usepackage[mathcal]{euscript}

\usepackage{cite}
\usepackage{hyperref}
\usepackage{enumitem}

\usepackage[ulem=normalem,draft]{changes}
%
%

%

\usepackage{xcolor}

\definecolor{viola}{rgb}{0.3,0,0.7}
\definecolor{ciclamino}{rgb}{0.5,0,0.5}
\definecolor{blu}{rgb}{0,0,0.7}
\definecolor{rosso}{rgb}{0.85,0,0}

\def\betti #1{{\color{magenta}#1}}
\def\bettib #1{{\color{magenta}#1}}
\def\juerg #1{{\color{green}#1}}
\def\pier #1{{\color{red}#1}}
\def\pcol #1{{\color{rosso}#1}}
\def\jupi #1{{\color{red}#1}} 
\def\revis #1{{\color{red}#1}}

\def\juerg #1{#1}
\def\betti #1{#1}
\def\bettib #1{{#1}}
\def\pier #1{{#1}}
\def\pcol #1{{#1}}
\def\jupi #1{#1}
\def\revis #1{#1}




\bibliographystyle{plain}


%
\newtheorem{theorem}{Theorem}[section]

\newtheorem{definition}[theorem]{Definition}

\finqui

\def\Beq{\Begin{equation}}
\def\Eeq{\End{equation}}

\def\Bthm{\Begin{theorem}}
\def\Ethm{\End{theorem}}

\def\Brem{\Begin{remark}\rm}
\def\Erem{\End{remark}}

\def\Bdim{\Begin{proof}}
\def\Edim{\End{proof}}
\def\Bcenter{\Begin{center}}
\def\Ecenter{\End{center}}
\let\non\nonumber




\def\step #1 \par{\medskip\noindent{\bf #1.}\quad}
\def\jstep #1: \par {\vspace{2mm}\noindent\underline{\sc #1 :}\par\nobreak\vspace{1mm}\noindent}

\def\Lip{Lip\-schitz}

\def\lhs{left-hand side}
\def\rhs{right-hand side}



\def\multibold #1{\def\arg{#1}%
  \ifx\arg\pto \let\next\relax
  \else
  \def\next{\expandafter
    \def\csname #1#1#1\endcsname{{\bf #1}}%
    \multibold}%
  \fi \next}

\def\pto{.}

\def\multical #1{\def\arg{#1}%
  \ifx\arg\pto \let\next\relax
  \else
  \def\next{\expandafter
    \def\csname cal#1\endcsname{{\cal #1}}%
    \multical}%
  \fi \next}


\def\multimathop #1 {\def\arg{#1}%
  \ifx\arg\pto \let\next\relax
  \else
  \def\next{\expandafter
    \def\csname #1\endcsname{\mathop{\rm #1}\nolimits}%
    \multimathop}%
  \fi \next}

\multibold
qwertyuiopasdfghjklzxcvbnmQWERTYUIOPASDFGHJKLZXCVBNM.

\multical
QWERTYUIOPASDFGHJKLZXCVBNM.

\multimathop
diag dist div dom mean meas sign supp .


\def\Accorpa #1#2 #3 {\gdef #1{\eqref{#2}--\eqref{#3}}%
  \wlog{}\wlog{\string #1 -> #2 - #3}\wlog{}}


\def\separa{\noalign{\allowbreak}}

\def\<#1>{\mathopen\langle #1\mathclose\rangle}
\def\norma #1{\mathopen \| #1\mathclose \|}

\def\I2 #1{\int_{Q_t}|{#1}|^2}
\def\IT2 #1{\int_{Q_t^T}|{#1}|^2}
\def\IO2 #1{\norma{{#1(t)}}^2}

\def\next{\\ & \quad}

\def\intQt{\int_{Q_t}}

\def\iO{\int_\Omega}

\def\dtt{\partial_{tt}}
\def\dt{\partial_t}
\def\dn{\partial_{\bf n}}

\def\checkmmode #1{\relax\ifmmode\hbox{#1}\else{#1}\fi}


\def\erre{{\mathbb{R}}}

\def\enne{{\mathbb{N}}}




\def\genspazio #1#2#3#4#5{#1^{#2}(#5,#4;#3)}
\def\spazio #1#2#3{\genspazio {#1}{#2}{#3}T0}

\def\L {\spazio L}
\def\H {\spazio H}
\def\W {\spazio W}

\def\C #1#2{C^{#1}([0,T];#2)}



\def\Lx #1{L^{#1}(\Omega)}
\def\Hx #1{H^{#1}(\Omega)}

\def\LQ #1{L^{#1}(Q)}

\def\Luno{\Lx 1}
\def\Ldue{\Lx 2}

\def\Huno{\Hx 1}
\def\Hdue{\Hx 2}


\def\LQ #1{L^{#1}(Q)}


\let\eps\varepsilon
\let\vp\varphi

\def\a{\alpha}	
\def\s{\sigma}  
\def\m{\mu}	    
\def\ph{\varphi}	

\def\h{\mathbbm{h}}

\let\TeXchi\chi                         
\newbox\chibox
\setbox0 \hbox{\mathsurround0pt $\TeXchi$}
\setbox\chibox \hbox{\raise\dp0 \box 0 }
\def\chi{\copy\chibox}



\def\mubar{\overline\mu}
\def\rhobar{\overline\rho}

\usepackage{amsmath}
\DeclareFontFamily{U}{mathc}{}
\DeclareFontShape{U}{mathc}{m}{it}%
{<->s*[1.03] mathc10}{}

\DeclareMathAlphabet{\mathscr}{U}{mathc}{m}{it}
\Begin{document}


%
\title{On the hyperbolic relaxation of the chemical potential in a  phase field  tumor growth model}

\author{}
\date{}
\maketitle
\Bcenter
\vskip-1.5cm
{\large\sc Pierluigi Colli$^{(1)}$}\\
{\normalsize e-mail: {\tt pierluigi.colli@unipv.it}}\\[0.25cm]
{\large\sc Elisabetta Rocca$^{(1)}$}\\
{\normalsize e-mail: {\tt elisabetta.rocca@unipv.it}}\\[0.25cm]
{\large\sc J\"urgen Sprekels$^{(2)}$}\\
{\normalsize e-mail: {\tt juergen.sprekels@wias-berlin.de}}\\[.5cm]
$^{(1)}$
{{\small Dipartimento di Matematica ``F. Casorati''}\\
{\small Universit\`a di Pavia}\\
{\small via Ferrata 5, I-27100 Pavia, Italy}\\[.3cm] 
$^{(2)}$
{\small Department of Mathematics}\\
{\small Humboldt-Universit\"at zu Berlin}\\
{\small Unter den Linden 6, D-10099 Berlin, Germany}\\[2mm]
{\small and}\\[2mm]
{\small Weierstrass Institute for Applied Analysis and Stochastics}\\
{\small Anton-Wilhelm-Amo-Strasse 39, D-10117 Berlin, Germany}\\[10mm]}
%
%
%
%
%
\Ecenter
\Begin{abstract}
\noindent 
In this paper, we study a phase field model for a tumor growth model of Cahn--Hilliard type in which  the often assumed parabolic relaxation of the chemical potential is replaced  by a hyperbolic one. \betti{We show that the resulting initial-boundary value problem is well posed and that its solutions depend continuously on two given functions: one appearing in the mass balance equation and one in the nutrient equation, representing, respectively, sources of drugs (e.g.~chemotherapy) and antiangiogenic therapy.} \pier{We also discuss regularity properties of the solutions. Moreover, in the case of a constant proliferation function, we rigorously analyze the asymptotic behavior as the coefficient of the inertial term tends to zero, establishing convergence to the corresponding viscous Cahn--Hilliard tumor growth model. Our results apply to a broad class of double-well potentials, including nonsmooth ones.}


\vskip3mm
\noindent {\bf Key words:}
Tumor growth models, singular potentials, hyperbolic relaxation, \pier{initial-boundary value problem, well-posedness, continuous dependence, regularity, asymptotic convergence}

\vskip3mm
\noindent {\bf AMS (MOS) Subject Classification:} \pier{%
35M33, 
37N25, 
35K57,  
35B30, 
35B65, 
35B40
}

\End{abstract}

%
\salta
\pagestyle{myheadings}
\newcommand\testopari{\sc Colli -- Rocca -- Sprekels}
\newcommand\testodispari{\sc Hyperbolic relaxation in a tumor growth model}
\markboth{\testopari}{\testodispari}
\finqui
%


\section{Introduction}
\label{INTRO}
\setcounter{equation}{0}

Let  $\a>0,~\pier{\tau}>0$, and let $\Omega\subset\erre^3$ denote some open and bounded domain having a smooth boundary $\Gamma=\partial\Omega$ with outward normal $\,{\bf n}\,$ and corresponding outward normal derivative $\dn$.  Moreover, we fix some final time $T>0$ and
introduce for every $t\in (0,T]$ the sets $Q_t:=\Omega\times (0,t)$ 
 and $\Sigma_t:=\Gamma\times (0,t)$,
 where we put, for the sake of brevity, $Q:=Q_T$ and $\Sigma:=\Sigma_T$.
We then consider the following initial-boundary vaue problem: 
\begin{align}
\label{ss1}
&\alpha\dtt\mu+\dt\ph-\Delta\mu=P(\ph)(\sigma+\chi(1-\ph)-\mu) - \h(\ph)u_1 &&\mbox{in }\,Q\,,\\[1mm]
\label{ss2}
&\pier{\tau}\dt\vp-\Delta\vp+F'(\vp)=\mu+\chi\,\sigma &&\mbox{in }\,Q\,,\\[1mm]
\label{ss3}
&\dt\sigma-\Delta\sigma=-\chi\Delta\vp-P(\ph)(\sigma+\chi(1-\ph)-\mu)+u_2 &&\mbox{in }\,Q\,,\\[1mm]
\label{ss4}
&\dn \mu=\dn\vp=\dn\sigma=0 &&\mbox{on }\,\Sigma\,,\\[1mm] 
\label{ss5}
&\mu(0)=\mu_0,\quad \dt\mu(0)=\mu_0',\quad\vp(0)=\vp_0,\quad \sigma(0)=\sigma_0 &&\mbox{in }\,\Omega\,.
\end{align}
\Accorpa\Statesys {ss1} {ss5}
The  system \Statesys\ constitutes a simplified and relaxed version of the four-species thermodynamically consistent model for tumor growth
originally proposed by Hawkins-Daruud et al.\ in \cite{HZO} that additionally includes chemotactic terms.
Let us briefly review the role of the occurring symbols. The primary variables $\ph, \m, \s$ denote the phase field, the associated chemical potential, and the nutrient concentration, respectively.
Furthermore, we stress that the additional term $\a\dtt\m$ 
is a hyperbolic regularization of equation \eqref{ss1},
whereas the term $\pier{\tau}\dt\ph$ is the viscosity contribution to the Cahn--Hilliard equation.
The key idea behind these regularizations originates from the fact that 
their presence allows us to take into account more general potentials $F$
whose derivatives $F'$, which play the role of a thermodynamic driving force in the model, may be singular and possibly nonregular.

The nonlinearity $P$ denotes a proliferation function, whereas the positive constant $\chi$
represents the chemotactic sensitivity. \pier{The terms containing $P(\varphi)$ in the system \eqref{ss1}--\eqref{ss5} model tumor cell proliferation. In \eqref{ss1}, 
the factor $P(\varphi)$ modulates the source of the chemical potential according to the local tumor density, coupling proliferation with nutrient availability ($\sigma$) and chemotaxis ($\chi(1-\varphi)$). In the third equation~\eqref{ss3}, 
the same term with minus sign accounts for nutrient consumption by proliferating cells. Thus, $P(\varphi)$ provides a natural coupling between the evolution of the phase field and the nutrient, reflecting the interplay between growth and resource uptake.}

\pier{About equation~\eqref{ss2}, as} is common in phase-field models, the function $F$
is assumed to have a double-well structure. Typical examples include the regular, logarithmic,
and double-obstacle potentials, which are respectively defined by
\begin{align}
\label{regpot}
&F_{\rm reg}(r) = \frac14 \left(1 - r^2\right)^2, \qquad r \in \mathbb{R}, \\[1mm]
\label{logpot}
&F_{\rm log}(r) = 
\begin{cases}
(1+r)\,\ln(1+r) + (1-r)\,\ln(1-r) - k_1 r^2, & r \in (-1,1),\\[1mm]
2\ln(2) - k_1, & r \in \{-1,1\},\\[1mm]
+\infty, & r \notin [-1,1],
\end{cases}
\\[1mm]
\label{2obspot}
&F_{\rm obs}(r) =
\begin{cases}
k_2 (1 - r^2), & r \in [-1,1],\\[1mm]
+\infty, & r \notin [-1,1],
\end{cases}
\end{align}
where $k_1 > 1$ and $k_2 > 0$, so that both $F_{\rm log}$ and $F_{\rm obs}$ are nonconvex.  
\pier{All these potentials can be written as the sum of a convex, lower-semicontinuous function (the main part) and a concave quadratic perturbation.}
Note that $F_{\rm log}$ is particularly relevant in applications, since $F'_{\rm log}(r)$ becomes
unbounded as $r \to \pm 1$. Moreover, in the case of \pier{$F_{\rm obs}$}, the second equation
\eqref{ss2} must be interpreted as a differential inclusion, where \pier{the derivative of the convex part of $F_{\rm obs}(\varphi)$ is understood in the sense of \juerg{subdifferentials}.}

In the above model equations, there are two functions that may 
serve as distributed controls acting in the phase and nutrient equations,
respectively. The control variable $u_1$, which is nonlinearly coupled to the state variable $\varphi$ in the phase equation \eqref{ss1}, models the application of a cytotoxic drug into 
the system; it is multiplied by a truncation function $\h(\cdot)$ in order to have the action
only in the spatial region where the tumor cells are located. For instance,
it can be assumed that $\h(-1)=0, \h(1)=1 ,\h(\ph)$ is in between if $-1<\ph<1$;
see \cite{GLSS, GL1, HKNZ, KL} for some insights on possible choices of $\h$.
On the other hand, the control $u_2$ can model  the supply of antiangiogenic therapies aiming at reducing the tumor vascularization 
\betti{(cf.~\cite{CGRS3} and \cite{CGLMRR1, CGLMRR2} for 
similar control terms in models of prostate tumor growth).}

\betti{Let us briefly recall the results already present in the literature on this class of models, which has been first introduced in \cite{HZO} in case $\a=\pier{\tau}=0$. 
As far as well-posedness is concerned, the above model has already been deeply investigated in the case $\a=\pier{\tau}=\chi=0$ (cf.~\cite{CGH,CGRS1,CGRS2,CRW, FGR}). Moreover, many variants of this model were considered and similar results were proven, \juerg{see, for instance,}
\cite{CSS1, EG, GL1, GL2, GLSS}.} \pier{In fact, a large body of literature is devoted to diffuse-interface and Cahn--Hilliard-type models for tumor growth. Foundational modelling contributions are due to Cristini, Lowengrub, Wise and collaborators~\cite{WLFC,CLLW,CL}, of course  enlightening the numerical and asymptotic investigations in~\cite{HZO,HKNZ}. Rigorous analytical studies began with~\cite{FGR}, later extended to multi-species mixtures and interactions in~\cite{DFRGM,FLR,FLRS}, to Brinkman-type or Darcy-type couplings in~\cite{EG,ACGW}, and to nonlocal and degenerate settings in~\cite{FLS}. Several works addressed chemotaxis, active transport, mechanical effects, or additional biological mechanisms, establishing well-posedness and qualitative properties for a range of tumour-growth systems~\cite{CGH,GL1,GL2,GLSS,GLR,GLS,SS}. Fractional and viscous variants of Cahn--Hilliard tumor-growth models, together with asymptotics, vanishing viscosities and vanishing relaxation limits, have been examined in~\cite{CGRS1,CGRS2,CGS24,CGS2X}. Additional results on long-time behaviour and stability can be found in~\cite{CRW,Yayla}. Many of these analyses rely on convexity and compactness tools.}

\pier{\juerg{The optimal} control of tumor-growth phase-field systems constitutes another important research direction. Early works on boundary and distributed control appeared in~\cite{CGRS3,GLR}. Control strategies incorporating chemotaxis, active transport, variable mobilities, and Keller--Segel dynamics were developed in~\cite{KL,SW,EG, AS, CGSS}. Further contributions included the optimal control theory and advanced optimality conditions in~\cite{EK,EK_ADV}, as well as refined analyses of treatment-time optimization and related asymptotics in~\cite{S_b,S_a,SigTime,You}. Singular logarithmic and double-obstacle potentials 
\juerg{were} addressed in~\cite{S,S_DQ}, whereas sparse controls and second-order conditions
\juerg{were} studied in~\cite{ST,CFGS}. Well-posedness, regularity, and asymptotic behavior for models relevant to control applications and including chemotaxis \juerg{were} developed in~\cite{CSS1}. These results collectively provide a rigorous framework for the design and optimization of therapeutic strategies governed by \revis{diffuse-interface} tumor-growth models. The authors of this paper intend to undertake a detailed analysis of the distributed control problems associated with the two controls $u_1$ and $u_2$ in a subsequent work.}

\pier{Concerning the hyperbolic relaxation of the chemical potential in the viscous Cahn--Hilliard equation (uncoupled from the nutrient and without mass sources), we refer to the recent contributions~\cite{CS,CS-contr}, which inspired the present work. In~\cite{CS}, well-posedness, continuous dependence, and 
regularity results were established, along with an analysis of the asymptotic behavior as the relaxation parameter $\alpha$ tends to \revis{zero}. A related optimal control problem was studied in~\cite{CS-contr}.}

\pcol{We now briefly outline the contents of the present paper. In Section~\ref{STATE}, we first prove the existence of a (weak) solution to the system \eqref{ss1}--\eqref{ss5}, together with a continuous dependence result of the solutions on the controls $u_1$ and $u_2$: all this is precisely stated in Theorem~\ref{THM:WEAK}. Section~\ref{REGUL} is devoted to establishing regularity results, see Theorem~\ref{THM:REGU}, which lead to the existence of a strong solution to \eqref{ss1}--\eqref{ss5} in a very general framework for the potentials, covering all the cases in \eqref{regpot}--\eqref{2obspot}. Finally, Section~\ref{Asy} addresses the asymptotic limit as $\alpha \searrow 0$ in the particular --- but still relevant --- case when the proliferation function $P(\varphi)$ is constant, providing a detailed proof of convergence to the corresponding system with $\alpha = 0$ (see Theorem~\ref{Convergence}). Moreover, we are able to prove an estimate of the difference of solutions in suitable norms with a precise rate of convergence (cf.~Theorem~\ref{Errest}).}

\section{General setting \pier{and well-posedness}}
\label{STATE}
\setcounter{equation}{0}
In this section, we introduce the general setting of our 
problem and state well-posedness results for the state system \eqref{ss1}--\eqref{ss5}. 
To begin with, for a Banach space $\,X\,$ we denote by $\|\cdot\|_X$
the norm in the space $X$ or in a power thereof, and by $\,X^*\,$ its dual space. 
The only exception from this rule applies to the norms of the
$\,L^p\,$ spaces and of their powers, which we often denote by $\|\cdot\|_p$, for
$\,1\le p\le \infty$. As usual, for Banach spaces $\,X\,$ and $\,Y\,$ that are contained in the same topological vector space, we introduce the linear space
$\,X\cap Y\,$ which becomes a Banach space when endowed with its natural norm $\,\|u\|_{X \cap Y}:=
\|u\|_X\,+\,\|u\|_Y$, for $\,u\in X\cap Y$.
Moreover, we introduce the spaces
\begin{align}
  & H := \Ldue \,, \quad  
  V := \Huno\,,   \quad
  W := \{v\in\Hdue: \ \dn v=0 \,\mbox{ on $\,\Gamma$}\}.
  \label{defHVW}
\end{align}
Furthermore, by $\,(\,\cdot\,,\,\cdot\,)$, $\,\Vert\,\cdot\,\Vert$ , and $\<\cdot,\cdot>$, we denote the standard inner product 
and related norm in $\,H$, as well as the dual product between $V$ and its dual $V^*$. We then have the dense and compact embeddings $\,V\subset H\subset  V^*$, with the standard identification  \,\,$\langle v,w \rangle =(v,w)\,\,$ for
 $v\in H$ and $w\in V$.  

\vspace{5mm}  
Throughout the paper, we make repeated use of H\"older's inequality, of the elementary Young's inequality
\begin{equation}
\label{Young}
a b\,\le \delta |a|^2+\frac 1{4\delta}|b|^2\quad\forall\,a,b\in\erre, \quad\forall\,\delta>0,
\end{equation}
as well as \juerg{of} the continuity of the embeddings $H^1(\Omega)\subset L^p(\Omega)$ for $1\le p\le 6$ and 
$\Hdue\subset C^0(\overline\Omega)$. Notice that the latter embedding is also compact, while this holds true
for the former embeddings only if $p<6$. We also introduce for $s\in(0,T]$ and \pcol{elements $w\in L^1(0,T;\Lx1)$ the notation 
\Beq
\label{convolution}
(1\ast w)(s):=\int_0^s w(\cdot,s')\,ds'\,,
\Eeq
that is, \pier{$(1\ast w)(s)\in \Lx1$} is the function that assigns to $x\in\Omega$ the value $\,\,\int_0^s
w(x,s')\,ds'$. Obviously, it holds that 
\begin{align}
\label{pieruwe}
&|(1\ast w)(s)|\,\le\,\revis{(1\ast |w|)(s)}
\quad\pier{\mbox{for a.e. $s\in (0,T)$}}\, .
\end{align}
Moreover, for $1\leq p,\, p'\leq \infty$ conjugate exponents and 
functions $v\in \pier{L^1(0,T;\Lx p)}$ and $w \in \pier{L^1(0,T;\Lx{p'})}$ we have that}
\begin{align}
\label{uwe}
\pcol{\bigg|\iO v(s)(1\ast w)(s)\bigg| \,\le\,
\pcol{\|v(s)\|_p \,\|w\|_{\pier{L^1(0,s:\Lx{p'})}}}\,
\quad\pier{\mbox{for a.e. $s\in (0,T)$}}\,.}
\end{align}

Finally, let us introduce a convention that will be tacitly employed throughout the paper:
the symbol  $C$ is used to indicate every constant
that depends only on the structural data of the problem (such as
$T$, $\Omega$, $\a$ or $\pier{\tau}$, the shape of the 
nonlinearities, and the norms of the involved functions),
so that its meaning may change from line to line.
If a parameter $\delta$ enters the computation, then
the symbol $C_\delta$ denotes constants
that additionally depend on~$\delta$.
On the contrary, precise constants that we refer 
to are denoted in a different way.

We now provide assumptions on the data of the problem. 
\begin{enumerate}[label={\bf (A\arabic{*})}, ref={\bf (A\arabic{*})}]
\item \label{const:weak}
	$\alpha,\pier{\tau} $ and $\chi$ are positive constants.
\item \label{F:weak}
	$F=F_1+F_2$ \pier{satisfies:} 
$\,F_1:\erre\to [0,+\infty]\,$ is convex \pier{and} lower semicontinuous with $\,F_1(0)=0$, while  $F_2 \in C^2(\erre)$ has a \Lip\ continuous derivative $F_2^{\,\prime}$. 
\item \label{P:weak}
	$P \in C^0(\erre)$ is nonnegative, bounded, and \Lip\ continuous.
\item \label{h:weak}
	$\h \in C^0(\erre)$ is nonnegative, bounded, and \Lip\ continuous.
\end{enumerate}
Let us note that all of the potentials  \eqref{regpot}--\eqref{2obspot} are 
admitted. 
In fact, the assumption \ref{F:weak} implies that the subdifferential $\partial F_1$ of $F_1$ is a 
maximal monotone graph in $\erre \times\erre$ with effective domain $D(\partial F_1 ) \subset D(F_1 ) $, and, since $F_1$
attains its minimum value~$0$ at $0$, it turns out that $0\in D(\partial F_1 )$ and $0\in\partial F_1(0)$. We also observe that the assumptions on $F_2$ imply that $F_2$ grows at most quadratically, that is, there are constants $\widehat c_1,\widehat c_2$ such that
\Beq
\label{quadratic}
|F_2(r)|\le \widehat c_1+\widehat c_2r^2 \quad\forall\,r\in\erre.
\Eeq  
\pier{Moreover,} we introduce the following notation: \pier{for} $r\in D(\partial F_1)$, we denote by $(\partial F_1)^\circ(r) $ the minimal section of $\partial F_1(r)$, that is,  the element of $\partial F_1 (r)$ having minimal modulus. \pier{Finally,
we extend the notations $F_1$, $\partial F_1$, $D(\partial F_1)$, and $(\partial F_1)^\circ$ to the corresponding functionals and operators induced on $L^2$ spaces.} 

Now, in the general setting  of~{\ref{const:weak}}--{\ref{h:weak}}, we are able 
to provide a well-posedness result for the system~\Statesys. 
First, we introduce our notion of a solution to \Statesys.

\begin{definition}
\label{DEF:WEAK}
A quadruple $(\m,\ph,\xi,\s)$ is called a solution to the initial-boundary value problem \Statesys\ if
\begin{align}
\mu&\in H^2(0,T;V^*)\cap W^{1,\infty}(0,T;H)\cap L^\infty(0,T;V),
\label{pier2-0}
\\
\ph & \in W^{1,\infty}(0,T;H)  \cap H^1(0,T;V) \cap L^\infty(0,T;W)\cap C^0(\overline Q), \label{pier2-1}
	\\
	\s & \in \H1 H \cap C^0([0,T];V) \cap \L2 {W}, \label{pier2-2}
	\\ 
	\xi & \in L^\infty(0,T;H), \label{pier2-3}
\end{align}
and if $(\m,\ph,\xi,\s)$ satisfies 
\begin{align}
	 \label{var:1}& \<\alpha\dtt\mu , v > + \iO\dt\ph \, v 
	+ \iO \nabla \mu \cdot \nabla v
	= \iO P(\ph)(\sigma+\chi(1-\ph)-\mu)v
	-\iO \h(\ph)u_1 v
	\nonumber\\
	& \qquad \hbox{for every $v \in V $ and a.e. in $(0,T)$,}\\[2mm]
	\label{var:2}&\pier{\tau}\dt\vp-\Delta\vp+\xi+F_2^{\,\prime}(\vp)=\mu+\chi\,\sigma, \quad \hbox{$\xi \in \partial F_1(\ph)$, \, a.e. in $\,Q$,}
 	\\
 	\label{var:3}
 	& \dt\sigma -\Delta\sigma
 	=-\chi\Delta\ph 
 	-  P(\ph)(\sigma+\chi(1-\ph)-\mu)
 	+ u_2	 \, \,\,\mbox{ a.e. in \,$Q$,}
\end{align}
as well as 
\begin{align}
	\label{var:4}
	\m(0)=\m_0, \quad
\dt\mu(0)=\mu_0',\quad
	\ph(0)=\ph_0, \quad
	\s(0)=\s_0, \quad \hbox{a.e. in $\Omega$}.
\end{align}
\end{definition}

It is worth noting that the homogeneous Neumann boundary conditions \eqref{ss4}
are considered in the \pier{conditions~\eqref{pier2-1} and \eqref{pier2-2} for $\ph$ and $\sigma$ (cf. the definition of the space $W$) and 
incorporated in the variational equality \eqref{var:1} for $\mu$,
when using the form $\iO \nabla \mu \cdot \nabla v  $.}  
\pier{Notice also that} the initial conditions~\eqref{var:4} make sense, since  \eqref{pier2-1} and \eqref{pier2-2} imply that $\ph,  \sigma \in
C^0([0,T];V)$, \pier{while, owing to \eqref{pier2-0}, it turns out that
$\mu \in   \C1{V^*}\cap \C0H$, and, consequently, $\dt \mu$ is at least weakly continuous from $[0,T]$ to $H$.}

We have the following result.
\begin{theorem}[Well-posedness]
\label{THM:WEAK}
Assume that \ref{const:weak}--\ref{h:weak} hold and \pier{let} the initial data  satisfy
\begin{align}
	\label{weak:initialdata}
	&\pier{{}\mu_0\in V, \quad \mu_0'\in H,} \quad
\s_0 \in V,  
	\quad
	\ph_0 \in W \pier{{}\cap D(\partial F_1)^\circ} 
\non\\
&\quad\hbox{\pier{with}} \quad 
	F_1(\ph_0) \in \Lx1, \quad \pier{(\partial F_1)^\circ} (\ph_0)\in H.
\end{align}
Moreover, suppose that 
\begin{align}
	\label{u:weak}
		(u_1, u_2) \in L^2(Q) \times L^2(Q).
\end{align}
Then there exists at least one solution $(\m,\ph,\xi,\s)$ in the sense of Definition \ref{DEF:WEAK}.
Moreover,~\pier{if 
\begin{align}
	\label{u:uniq}
		(u_1, u_2) \in \pier{\L2{\Lx3}} \times L^2(Q)
\end{align}
in addition to \eqref{u:weak},} then the  solution is unique.
Furthermore, let $(\m_i,\ph_i,\xi_i,\s_i)$, $i=1,2$, be two solutions to \Statesys\ associated 
with the \pier{data} $(u_1^i,u_2^i) \in \pier{\L2{\Lx3}} \times L^2(Q)$,
$i= 1,2 $. Then there exists a constant $K_1>0$, which depends only on the data of the system, such that
\begin{align}
	& \non
	\pier{ \|\m_1-\m_2\|_{\L\infty H}}
	+ \|1\ast(\m_1-\m_2)\|	_{\L \infty V}
		\\
	& \quad \non
	+ \norma{\ph_1-\ph_2}_{\L\infty H \cap \L2 V}
	+ \norma{\s_1-\s_2}_{\L\infty H \cap \L2 V}
	\\
&
\label{cont:dep:weak}
	\leq
	K_1  \Big(
	 \norma{u_1^1-u_1^2}_{\L2 H}
	+  \norma{u_2^1-u_2^2}_{\L2 H}
	\Big).
\end{align} 
\end{theorem}
Before entering the proof, let us remark that the above result is very general and includes also the cases of singular and 
nonsmooth potentials such as the double obstacle potential defined in \eqref{2obspot}.
\pier{We also note that the assumption $F_1(\ph_0) \in \Lx1$ stated in \eqref{weak:initialdata} actually follows from other requirements on $\ph_0$ thanks to the subdifferential property
$$ 
\iO F_1(\ph_0) \leq \iO F_1(v) + ( \pier{(\partial F_1)^\circ} (\ph_0), \ph_0 -v) \quad \hbox{for every $v\in H$}.
$$
About} the explicit dependencies of the constant $K_1$ \pier{in ~\eqref{cont:dep:weak}}, we invite the reader to follow the proof of 
the estimate \pier{given} below. 
\begin{proof}  
The existence proof is rather standard, since similar  arguments have already been used in previous contributions. 
Hence, for that part, we proceed rather formally, just employing 
the Yosida approximation of $\partial F_1$ for our estimates without recurring to finite-dimensional
approximation techniques like the Faedo--Galerkin scheme.
Hence, we introduce the Yosida regularization of $ \partial F_1$. 
For $\eps>0$, let $F_{1,\eps}$ denote the Moreau--Yosida
approximation of $F_1$ at the level $\eps$. 
It is well known (see, e.g., \cite{Brezis}) 
that the following conditions are satisfied:
\begin{align}
&0\leq F_{1,\eps}(r) \le F_1(r) \quad\mbox{for all }\,r\in\erre\pier{\,;}
\label{pier3} \\
&F^{\,\prime}_{1,\eps} \mbox{ is Lipschitz continuous on $\erre$} \non\\
&\quad \mbox{with \Lip\ constant \pier{$1/\eps$}, and $F^{\,\prime}_{1,\eps}(0)=0$\pier{\,;}}
\label{pier4}\\
&|F^{\,\prime}_{1,\eps}(r)|\le |(\partial F_1)^\circ (r)| \,\,\mbox{ and }\,\,\lim_{\eps\searrow 0}\,F^{\,\prime}_{1,\eps}(r)=(\partial F_1)^\circ (r)
\quad\mbox{for all }\,r\in D(\partial F_1)\, . 
 \label{pier5} 
\end{align} 
We now study the \revis{$\varepsilon$-approximating} problem that results from the system \eqref{var:1}--\eqref{var:4} if  
$(\partial F_1) (r)$ is replaced by $F^{\,\prime}_{1,\eps}$ and the inclusion in  \eqref{var:2} reduces to an equality. Namely, we argue on
\begin{equation}
\label{pier6}
\pier{\tau}\dt\vp-\Delta\vp+F^{\,\prime}_{1,\eps} (\vp) +F_2^{\,\prime}(\vp)=\mu+\chi\,\sigma \quad \hbox{a.e. in $Q$.}
\end{equation}
The existence of a solution to the \revis{$\varepsilon$-approximating} problem thus obtained can be shown by means of a Faedo--Galerkin approximation using the eigenvalues 
$\{\lambda_j\}_{j\in\enne}$ and eigenfunctions $\{e_j\}_{j\in\enne}$ of the  eigenvalue problem  $\,\,-\Delta e_j=\lambda_je_j\,\,$ in $\Omega$, \,$\dn e_j=0\,$ on $\Gamma$. In order not to overload the exposition,  we avoid here to write the Faedo--Galerkin system explicitly. Instead, we just provide the relevant a priori estimates
by performing the estimations directly on the solution to the $\varepsilon$-approximating system.  Notice that the following estimates, while being only formal for the $\eps$-approximating system, are fully justified on the level of the Faedo--Galerkin approximations. 
For the sake of simplicity, we still denote by $(\m,\ph,\xi,\s)$, with $\xi = F^{\,\prime}_{1,\eps} (\vp)$,
the solution to the \revis{$\varepsilon$-approximating} system in place of $(\m_\eps,\ph_\eps,\xi_\eps,\s_\eps)$; the correct notation will be reintroduced at the end of each estimate. Before entering the estimates, we note that it follows from \eqref{weak:initialdata} that $\ph_0\in C^0(\overline \Omega)$, and  we \pier{conclude} from the assumption \ref{F:weak} that 
$F_2(\ph_0)\in\Luno$ and $F_2^{\,\prime}(\ph_0)\in H$. 

\noindent
{\sc First Estimate:}
Let $t\in(0,T]$ be arbitrary. First, we test \eqref{var:1} by $\dt\mu$
\pier{and integrate over $[0,t]$.} 
Using \ref{P:weak}, \ref{h:weak}, \eqref{weak:initialdata} and Young's inequality, we then obtain that 
\begin{align}
\label{esti1}
&\frac\a 2\|\dt\mu(t)\|^2 \pier{{}+\int_{Q_t}\dt\ph\dt\mu 
+ \frac 12\iO |\nabla \mu(t)|^2 } \non\\
&=\frac \a 2\|\mu_0'\|^2 + \pier{\frac 12\iO |\nabla \mu_0|^2 }
+\int_{Q_t}\bigl[P(\ph)(\s+\jupi{\chi (1-\ph)}-\mu) -\h(\ph)u_1\bigr]\dt\mu
\,,
\end{align}
\pier{\juerg{where the last integral} can be estimated as follows: 
\begin{align}
\label{esti1pier}
&\int_{Q_t}\bigl[P(\ph)(\s+\jupi{\chi (1-\ph)}-\mu)-\h(\ph)u_1\bigr]\dt\mu
\non\\
&\le C\intQt \bigl(|\s|^2 \jupi{{}+1{}} +|\ph|^2+ |\mu|^2+ |u_1|^2+|\dt\mu|^2\bigr)\,.
\end{align}
Next,} we differentiate \eqref{pier6} with respect to $t$, test by $\dt\ph$, and integrate over $Q_t$. We then find the identity
\begin{align}
\label{esti2}
&\frac{\pier{\tau}} 2\|\dt\ph(t)\|^2+\intQt |\nabla\dt\ph|^2+\intQt F^{\,\prime\prime}_{1,\varepsilon}(\ph)|\dt\ph|^2\non\\
&=\frac{\pier{\tau}} 2\|\dt\ph(0)\|^2-\intQt F^{\,\prime\prime}_2(\ph)|\dt\ph|^2 +\pier{\intQt\dt\mu \dt\ph}+\chi\intQt \dt\s\dt\ph \,.
\end{align}
Now, writing \eqref{pier6} at $t=0$, we see that
$$
\dt\ph(0)=\pier{\tau}^{-1}\bigl(\Delta\ph_0 -F^{\,\prime}_{1,\varepsilon}(\vp_0)-F_2^{\,\prime}(\vp_0)+\mu_0+\chi\s_0\bigr),
$$
and we infer from \pier{\eqref{weak:initialdata}, \eqref{pier5} and \ref{F:weak} that
$\|\dt\ph(0)\|$ is uniformly bounded}. Moreover, we have $\,F_{1,\varepsilon}^{\,\prime\prime}(\ph)\ge 0\,$ so that the last term on the \lhs\ is nonnegative, while $F_2^{\,\prime\prime}(\ph)$ is bounded. Hence,
applying Young's inequality to the \pier{fourth} term on the \rhs\ and then adding the inequalities \eqref{esti1} and \eqref{esti2}, thus cancelling the terms involving
$\dt\ph\dt\mu$, we arrive at
\begin{align}
\label{esti3}
&\frac\a 2\|\dt\mu(t)\|^2 + \pier{\frac 12\iO |\nabla \mu(t)|^2} +\frac\tau 2\|\dt\ph(t)\|^2
+\intQt |\nabla\dt\ph|^2
\non\\
&\le \,C+C\intQt \bigl(|\mu|^2 + |\ph|^2 + |\s|^2 + |\dt\mu|^2 + |\dt\ph|^2
\bigr) +\pier{\frac 1{16}}\intQt |\dt\s|^2\,. 
\end{align}

Next, we test \eqref{var:3} by $\dt\s$, integrate over $(0,t)$, and \pier{add to}
both sides the expression  $\,\,\frac 12\|\s(t)\|^2-\frac 12 \|\s_0\|^2=\intQt
\s\dt\s$. We then find that 
\begin{align}
&\intQt|\dt\s|^2 + \frac 12 \|\s(t)\|_V^2
= \frac 12\|\s_0\|_V^2+\chi\intQt\nabla
\ph\cdot\nabla\dt\s \non\\
&\quad{}+\intQt\bigl[-P(\ph)(\s+\jupi{\chi (1-\ph)}-\mu)+\s+u_2\bigr]\dt\s\,. \non
\end{align}
Now observe that, \pier{integrating by parts in time and applying Young's inequality,}
\begin{align}
&\chi\intQt\nabla\ph\cdot\nabla\dt\s =\chi\iO\nabla\ph(t)\cdot\nabla\s(t)
-\chi\iO\nabla\ph_0\cdot\nabla\s_0-\chi\intQt\nabla\dt\ph\cdot\nabla\s
\non\\
&\le\,C \pier{{}+\frac 14\iO |\nabla\s(t)|^2  +\chi^2 \iO|\nabla\ph(t)|^2}
+\frac 12\intQt|\nabla\dt\ph|^2+\frac{\chi^2}2\intQt|\nabla\s|^2\,.\non
\end{align}
Hence, also applying Young's inequality to the last term on the \rhs\ of the 
penultimate identity, we can infer that
\begin{align}
\label{esti4}
&\frac12 \intQt|\dt\s|^2 + \frac 14 \|\s(t)\|_V^2\le C+\chi^2\iO|\nabla\ph(t)|^2\non\\
&\quad{}+\frac 12\intQt|\nabla\dt\ph|^2
+C\intQt\pier{\bigl(|\nabla\s|^2 + |\s|^2 \jupi{{}+1{}} +|\ph|^2+ |\mu|^2+ |u_2|^2 \bigr)}
\,.
\end{align}
Finally, we test \eqref{pier6} by $\dt\ph$, integrate over $Q_t$, and add to both 
sides the expression $\,\,\frac 12 \|\ph(t)\|^2-\frac 12\|\ph_0\|^2=\intQt \ph\dt\ph$. 
\pier{In view of \ref{F:weak} and \eqref{pier3}, and taking the quadratic growth of 
$F_2$ into account, we obtain that}
\begin{align}
\label{esti5}
&\pier{\pier{\tau} \intQt |\dt\ph|^2} +\frac 12\|\ph(t)\|_V^2+\pier{\iO F_{1,\varepsilon}(\ph(t))} 
\non\\
&=\frac 12\|\ph_0\|_V^2
+\iO \bigl( F_{1,\varepsilon}(\ph_0) + F_2(\ph_0)\bigr) 
-\iO F_2(\ph(t)) +\intQt (\mu+\chi\s+\ph)\dt\ph\non\\
&\le \pier{C
+\iO F(\ph_0)+C\left(1+\|\ph_0\|^2+\intQt\ph\dt\ph\right)}+\intQt (\mu+\chi\s+\ph)\dt\ph\non\\
&\le\,C+C\intQt\bigl(|\mu|^2+|\ph|^2+|\s|^2+|\dt\ph|^2\bigr)\,,
\end{align}
\pier{where we have used Young's inequality and \eqref{weak:initialdata} as well. Note moreover that the third term on first line of \eqref{esti5} is nonnegative (cf.~\eqref{pier3}).}

At this point, we multiply \eqref{esti5} by $\,4\chi^2\,$ and add the resulting inequality  to the
sum of the inequalities  \eqref{esti3} and \eqref{esti4}. \pier{We infer that 
\begin{align}
\label{esti6pier}
&\frac\a 2\|\dt\mu(t)\|^2 + \pier{\frac 12\iO |\nabla \mu(t)|^2} +\frac\tau 2\|\dt\ph(t)\|^2 +\frac{\chi^2}2\|\ph(t)\|_V^2 \non \\
&\quad +  4 \tau \chi^2 \intQt |\dt\ph|^2 + \frac12 \intQt |\nabla\dt\ph|^2 
+ \frac14 \intQt|\dt\s|^2 + \frac 14 \|\s(t)\|_V^2
\non\\
&\le \,C + C\intQt \bigl(|\mu|^2 + |\ph|^2 + |\s|^2 + |\nabla\s|^2 + |\dt\mu|^2 + |\dt\ph|^2
\bigr) \,. 
\end{align}
Now, note that the term $ C\intQt |\mu|^2$ on the \rhs\ can be estimated using the identity 
$\mu (s) =\mu_0 + \int_0^s \dt \mu $, $s\in (0,t)$,  so that 
$ C\intQt |\mu|^2 \leq C + C\intQt  |\dt\mu|^2$. It is then easily seen that 
the  inequality thus obtained from \eqref{esti6pier}
admits the application of Gronwall's lemma,
and we finally can deduce that}
\begin{align}
\label{smokeyuno}
&\|\mu_\varepsilon\|_{W^{1,\infty}(0,T;H)\cap L^\infty(0,T;V)} \,+\,
\|\ph_\varepsilon\|_{W^{1,\infty}(0,T;H)\cap H^1(0,T;V)}\non\\
&\quad + \|\s_\varepsilon\|_{H^1(0,T;H)\cap L^\infty(0,T;V)}\,\le\,C\,.
\end{align}

\noindent
{\sc Second Estimate:}
Next, we multiply \eqref{pier6}, written for a fixed $t\in (0,T]$, by $\,\,-\Delta\ph(t)\,\,$ and integrate over $\Omega$. This yields the identity
\Beq
\nonumber
\|\Delta\ph(t)\|^2+\iO F_{1,\eps}^{\,\prime\prime}(\ph(t))|\nabla\ph(t)|^2
=\iO\pier{(\pier{\tau}\dt\ph + F_2^{\,\prime}(\ph) - \chi\s -\mu)}(t)\Delta\ph(t)
\Eeq
and, using Young's inequality, \pier{the monotonicity of $F_{1,\eps}^{\,\prime}$, 
the bound in \eqref{smokeyuno} which implies \juerg{an} $\L\infty H$-bound for $(\pier{\tau}\dt\ph + F_2^{\,\prime}(\ph) - \chi\s -\mu)$,  and the elliptic regularity theory, we plainly} deduce that
\Beq
\label{smokeydue}
\|\ph_\eps\|_{L^\infty(0,T;W)} \,\le\,C\,.
\Eeq
But then, by comparison in equation \eqref{pier6}, \pier{we also realize that}
\Beq
\label{smokeytre}
\|F_{1,\eps}'(\ph_\eps)\|_{L^\infty(0,T;H)}\,\le\,C\,,
\Eeq
\pier{while comparison in \eqref{var:3}, along with \eqref{smokeyuno} and elliptic regularity again, yields that}
\Beq
\|\sigma_{\varepsilon}\|_{L^2(0,T;W)}\,\le\,C\,.
\label{smokey3bis}
\Eeq
Finally, \pier{it follows from a comparison of terms} in \eqref{var:1} that
\Beq
\label{smokeyfour}
\|\betti{\mu}_{\varepsilon}\|_{H^2(0,T;V^*)}\,\le\,C\,.
\Eeq

\noindent 
{\sc Passage to the limit as $\eps\searrow0$:}
Now, let for every $\eps>0$ the triple $(\m_\eps,\ph_\eps,\s_\eps)$ be a solution to the problem \eqref{var:1}, \eqref{pier6}, \eqref{var:3}, \eqref{var:4} with the regularity \eqref{pier2-0}--\eqref{pier2-2}. Observe that the constants $C$ occurring in the proof of the estimates \eqref{smokeyuno}--\eqref{smokeyfour} are all independent of $\eps$. Hence it follows from standard weak and weak-star compactness results that there are functions 
$\m,\ph,\s,\xi $  such that  
\begin{align}
\label{pier10-1}
& \mu_\eps \to \mu \quad \hbox{weakly star in } \,  H^2(0,T;V^*) \cap 
W^{1,\infty}(0,T;H) \cap L^\infty(0,T;V),   
\\
\label{pier10-2}
& \ph_\eps \to \ph \quad \hbox{weakly star in } \,   W^{1,\infty}(0,T;H) \cap H^1(0,T;V) \cap L^\infty(0,T;W), 
\\
\label{pier10-3}
&\s_\eps \to \s \quad \hbox{weakly star in } \,   H^1(0,T;H) \cap L^\infty(0,T;V) \cap 
L^2(0,T;W),  
\\
\label{pier10-4}
& F^{\,\prime}_{1,\eps}(\ph_\eps) \to \xi \quad \hbox{weakly star in } \,   L^\infty(0,T;H), 
\end{align}
as $\eps \searrow 0$, possibly along a subsequence. By virtue of \eqref{pier10-2} and 
the Aubin--Lions--Simon lemma (see, e.g., \cite[Sect.~8, Cor.~4]{Simon}), \pier{as $W$ is compactly embedded into $C^0(\overline \Omega)$,} we deduce that  
\Beq
\ph_\eps\to \ph \,\mbox{ strongly in \,$C^0(\overline Q)$}, \non
\Eeq
 whence, by
Lipschitz continuity, also
\begin{equation*}
P(\ph_\eps)\to P(\ph), \quad \h(\ph_\eps)\to\h(\ph),
\quad F_2^{\,\prime}(\ph_\eps)\to F_2^{\,\prime}(\ph), \quad\mbox{all  strongly in \,$C(\overline Q)$.}
\end{equation*}
\pier{On the other hand, by the same tool, we have that
\begin{align}
\non 
& \mu_\eps \to \mu \quad \hbox{strongly in } \, \C1{V^*}\cap \C0 H ,  
\\
\non
&\s_\eps \to \s \quad \hbox{strongly in } \,   \C0H \cap 
L^2(0,T;V),
\end{align}
as consequences of \eqref{pier10-1} and \eqref{pier10-3}.}

We then can pass to the limit in the \pier{respective variational equality~\eqref{var:1} and equation~\eqref{var:3}, in particular, for the product term $P(\ph_\eps)(\sigma_\eps+\chi(1-\ph_\eps)-\mu_\eps)$.}
This is also possible in \eqref{pier6} in order to obtain \pier{\eqref{var:2} \juerg{in} the limit:} indeed,  
the inclusion in~\eqref{var:2} 
results as a consequence of \eqref{pier10-4} and the
maximal monotonicity of $\partial F_1$, since we can apply, e.g., 
\cite[\pier{Prop.~2.2,} p.~38]{Barbu}. Finally, the initial conditions 
\eqref{var:4} can be readily obtained by observing that \pier{we have at least strong convergence in $C^0([0,T]; V^*)$ for all 
\juerg{of the variables}
$\mu_\eps, \dt\mu_\eps, \ph_\eps, \sigma_\eps$}. With this, the existence part of the proof is complete.

\vspace{2mm} \noindent
{\sc Uniqueness and Continuous Dependence:} 
We now show the continuous dependence estimate \eqref{cont:dep:weak}, which implies the uniqueness of the solution, in particular. To this end, let $(u_1^i,u_2^i)\in  \revis{\L2{\Lx3} \times L^2(Q)}$, $i=1,2$, be given, and let  $(\mu_i,\ph_i,\xi_i,\s_i)$, $i=1,2$, be associated solutions in the sense of Definition 2.1. We then introduce the abbreviating notation
\begin{align}
&\mu:=\mu_1-\mu_2, \quad \ph:=\ph_1-\ph_2,\quad \xi:=\xi_1-\xi_2,
\quad \s:=\s_1-\s_2,
\non\\
&u_1:=u_1^1-u_1^2, \quad u_2:=u_2^1-u_2^2. \non
\end{align}
We then see that the differences $(\mu,\ph,\xi,\s)$ satisfy the identities
\begin{align}
 \label{diff1}
& \langle \alpha\dtt\mu , v \rangle + \iO\dt\ph v 
	+ \iO \nabla \mu \cdot \nabla v
	= \iO (P(\ph_1)-P(\ph_2))(\sigma_1+\jupi{\chi (1-\ph_1)}-\mu_1)v
\non\\
&\quad{} +\iO P(\ph_2)(\s \jupi{{}- \chi\, \ph -{}} \m)v
-\iO (\h(\ph_1)-\h(\ph_2))u_1^1 v -\iO\h(\ph_2)u_1 v\nonumber\\
	& \qquad{} \hbox{for every $v \in V $ and a.e. in $(0,T)$,}\\[2mm]
\label{diff2}
&\pier{\tau}\dt\vp \revis{{}-\Delta{}}
\vp+\xi+F_2^{\,\prime}(\vp_1)-F_2^{\,\prime}(\ph_2)=\mu+\chi\,\sigma, 
\non\\
&\qquad{} \hbox{$\xi_i \in \partial F_1(\ph_i)$, \ $i=1,2$, \ a.e. in $\,Q$,}
 	\\[2mm]
\label{diff3}
 	& \dt\sigma -\Delta\sigma
 	=-\chi\Delta\ph 
 	-  (P(\ph_1)-P(\ph_2))(\sigma_1+\jupi{\chi (1-\ph_1)}-\mu_1)
\non\\
&\hspace*{22mm} - P(\ph_2)(\s\jupi{{}- \chi\, \ph -{}}\mu)
 	+ u_2	 \, \,\,\mbox{ a.e. in \,$Q$,}
\\[2mm]
\label{diff4}
&\mu(0)=0,\quad \dt\m(0)=0,\quad \ph(0)=0,\quad \s(0)=0,\quad
\mbox{ a.e. in $\,\Omega$}. 
\end{align}
Now recall that $\ph_i\in C^0(\overline Q)$, $i=1,2$. Hence, there is some constant 
$L>0$, which only depends on $R:= \max\,\{\|\ph_1\|_{C^0(\overline Q)}\,,\,
\|\ph_2\|_{C^0(\overline Q)}\}$, such that
\Beq
\label{differ1}
|P(\ph_1)-P(\ph_2)| + |\h(\ph_1)-\h(\ph_2)| + |F_2^{\,\prime}(\ph_1)-F_2^{\,\prime}(\ph_2)| \,\le \,L|\ph| \quad
\mbox{a.e. in \,$Q$}.
\Eeq
Let $t\in (0,T]$ be arbitrary. We multiply \eqref{diff2} by $\ph$ and integrate over $Q_t$ to obtain from \eqref{differ1} and Young's inequality that
\begin{align}
\label{differ2}
&\frac{\pier{\tau}} 2\|\ph(t)\|^2+\intQt|\nabla\ph|^2+\intQt \xi \ph\non\\
&=\intQt \bigl(\chi\s-(F_2^{\,\prime}(\ph_1)-F_2^{\,\prime}(\ph_2))\bigr)\ph+\intQt \mu\ph \non\\
&\le C\intQt (|\s|^2 + |\ph|^2)+\intQt \mu\ph\,.
\end{align}
\pier{We} observe that the last term on the \lhs\ \pier{of \eqref{differ2} is 
nonnegative on account of} the monotonicity of~$\partial F_1$.

\pier{Next,} let $M>0$ denote a constant that will be specified later. We multiply \eqref{diff3} by $M\s$
and integrate over $Q_t$. Owing to the boundedness of $P$, and invoking \eqref{differ1}, we then obtain the estimate 
\begin{align}
\label{uff0}
&\frac M 2 \|\s(t)\|^2 + M\intQt |\nabla\s|^2\non\\
&\le\,M\chi\intQt \nabla\ph\cdot\nabla \s \,+\,ML\intQt|\ph|\,\bigl(\jupi{1+{}}|\mu_1| + |\ph_1| + |\s_1|\bigr)
|\s| \non\\
&\quad +MC\intQt \bigl(|\mu|^2 + |\ph|^2+|\s|^2\bigr)+M\intQt |u_2|^2\,,
\end{align}
where \pier{we} have applied Young's inequality to the last two terms on the \rhs.
\pier{Besides,} we observe that
\Beq
\label{uff1}
M\chi\intQt \nabla\ph\cdot\nabla \s\,\le\,\frac M2\intQt|\nabla\s|^2+\frac{M\chi^2}2
\intQt|\nabla\ph|^2\,.
\Eeq
Finally, we recall that, \pier{thanks} to the continuity of the embedding $V\subset L^4(\Omega)$ and
the regularity \pier{properties}~\eqref{pier2-0}--\eqref{pier2-2}, the functions $\mu_1,\ph_1,\s_1$ are all
bounded in $L^\infty(0,T;L^4(\Omega))$. Therefore, we can infer that the second term on the \rhs\ of \eqref{uff0}, which we denote by $I$, can be estimated as follows: 
\begin{align}
\label{uff2}
|I|\,&\le\, MC\int_0^t \|\ph(s)\|_4 \,\bigl(\jupi{1+{}}\|\mu_1(s)\|_4 + \|\ph_1(s)\|_4 +\|\s_1(s)\|_4\bigr)\,
\|\s(s)\|\,ds\non\\
{}&\le\, MC\int_0^t \|\ph(s)\|_V \,\|\s(s)\|\,ds\non\\
{}&\le\,MC\intQt|\ph|^2 \,+\,\frac{M\chi^2}2\intQt |\nabla\ph|^2\,+MC\intQt|\s|^2\,.
\end{align}
Hence, combining \eqref{uff0}--\eqref{uff2}, we have shown the estimate
\begin{align}
\label{differ3}
&\frac M2\|\s(t)\|^2+\frac M2\intQt|\nabla\s|^2\non\\
&{}\le\,
M\chi^2\intQt|\nabla\ph|^2 + MC\intQt \bigl(|\mu|^2+|\ph|^2+|\s|^2+|u_2|^2\bigr)\,.
\end{align}
It remains to treat the identity \eqref{diff1} which we integrate with respect to time over \pier{$[0,s]$} for
$s\in (0,t]$. Then we insert $v=\mu(s)$ and integrate over $[0,t]$, arriving at the identity
\begin{align}
\label{uff3}
&\frac \alpha 2\|\mu(t)\|^2 + \frac 12 \pier{\iO |\nabla(1\ast \mu)\pier{(t)}|^2}
+\intQt \mu\ph\non\\
&=\,\intQt \bigl[1\ast \bigl((P(\ph_1)-P(\ph_2))(\s_1+\jupi{\chi (1-\ph_1)}-\mu_1)\bigr)\bigr]\,\mu\non\\
&\quad + \intQt \bigl[1\ast \bigl( P(\ph_2)(\s \jupi{{}-\chi \ph}-\mu)\bigr)\bigr]\,\mu 
-\intQt \bigl[1\ast \bigl((\h(\ph_1)-\h(\ph_2))\,u_1^1\bigr)\bigr]\,\mu\non\\
&\quad -\intQt \bigl[1\ast \bigl(\h(\ph_2)u_1\bigr)\bigr]\,\mu \,\,\,=: I_1+I_2+I_3+I_4\,,
\end{align}
with obvious meaning, where we have used the notation introduced in \eqref{convolution}. 
We estimate the terms on the \rhs\ of \eqref{uff3} individually, where we make repeated use of 
the H\"older and Young inequalities,  the \pcol{estimates \eqref{pieruwe} and \eqref{uwe} with $p=2$, and  the continuous embedding $V\subset \pier{L^q(\Omega)} $ \pier{for $1\leq q \leq 6$}}.   
At first, \pier{by the H\"older and Young inequalites, as} $u_1^1\in\pier{\L2{\Lx3}}$ it is readily seen that
\begin{align}
\label{uff-uff1}
I_3\,& \le\,L \int_0^t \|\mu(s)\|\,\pier{\|\, |\ph|\, |u_1^1|\, \|_{L^1(0,s;H)}}\,ds
\non\\
&\pier{{}\le \, 
L \int_0^t \|\mu(s)\|\int_0^s \|u_1^1(s')\|_{3} \|\ph(s')\|_{6} \,ds'\, ds}
\non\\
&\pier{{}\le \, 
C \|u_1^1 \|_{L^2(0,T;\Lx3)} \int_0^t \|\mu(s)\| \,  \|\ph\|_{L^2(0,s;V)} \, ds}
\non\\
&\pier{{}\le\,\frac14 \intQt|\nabla\ph|^2+\frac14\intQt|\ph|^2+C\int_0^t \|\mu(s)\|^2 \, ds\,.}
\end{align}
Similar reasoning, using the boundedness of both $P$ and $\h$, yields
\begin{align}
\label{uff-uff2}
I_2+I_4\,&\pier{{}\le \,C \int_0^t \|\mu(s)\|\,\pier{\|\, |\sigma|+|\ph|+|\mu| + |u_1|\, \|_{L^1(0,s;H)}}\,ds}
\non\\
&\pier{{}\le \,C \int_0^t \|\mu(s)\|\,\pier{\|\, |\sigma|+|\ph|+|\mu| + |u_1|\, \|_{L^2(0,t;H)}}\,ds}
\non\\
&\le\,C\intQt \bigl(|\mu|^2+|\s|^2+|\ph|^2+|u_1|^2\bigr)\,.
\end{align}
Finally, using \eqref{uwe}, \eqref{differ1}, as well as H\"older's and Young's inequalit\pier{ies}, and invoking the fact that $\s_1,\ph_1,\mu_1\in L^\infty(0,T;L^4(\Omega))$,
we find \juerg{the chain} of inequalities
\begin{align}
\label{uff-uff3}
I_1\,&\pier{{}\le \,C \int_0^t \|\mu(s)\|\,\|\,|\ph|( |\sigma_1|\jupi{{}+1{}}+|\ph_1|+|\mu_1|) \|_{L^1(0,s;H)} \,ds}
\non\\
&\pier{{}\le \, 
C \int_0^t \|\mu(s)\|\int_0^s \|(|\sigma_1|\jupi{{}+1{}}+|\ph_1|+|\mu_1|)(s')\|_{4} \|\ph(s')\|_{4} \,ds'\, ds}
\non\\
&\pier{{}\le \, 
C\,  \|\,\jupi{{}1+{}}\,|\mu_1| + |\ph_1| +|\s_1|\, \|_{L^\infty(0,T;\Lx4)} \int_0^t \|\mu(s)\| \,  \|\ph\|_{L^2(0,s;V)} \, ds}
\non\\
&\pier{{}\le\,\frac14 \intQt|\nabla\ph|^2+\frac14\intQt|\ph|^2+C\int_0^t \|\mu(s)\|^2 \, ds\,.}
\end{align}
Combining the estimates \eqref{uff3}--\eqref{uff-uff3}, we thus have shown the estimate
\begin{align}
\label{differ4}
&\pier{{}\frac \alpha 2\|\mu(t)\|^2 + \frac 12 \iO |\nabla(1\ast \mu)\pier{(t)}|^2 +\intQt \mu\ph}
\non\\
&\pier{{}\le\,\frac12 \intQt|\nabla\ph|^2 +C \intQt\bigl(|\mu|^2+|\s|^2+|\ph|^2+|u_1|^2\bigr).}
\end{align}
At this point, we add the inequalities \eqref{differ2}, \eqref{differ3} and \eqref{differ4}, obtaining the estimate
\begin{align}
\label{differ5}
&\frac{\pier{\tau}}2\|\ph(t)\|^2+\frac M2 \|\s(t)\|^2 + \pier{{}\frac \alpha 2\|\mu(t)\|^2 + \frac 12 \iO |\nabla(1\ast \mu)\pier{(t)}|^2}
\non\\
&\quad\pier{{}+\frac12 \bigl(1-2M\chi^2\bigr)}\intQt|\nabla\ph|^2+\frac M2\intQt|\nabla\s|^2\non\\
&\le\,\pier{C_{M}}\intQt\bigl(|\mu|^2+|\s|^2+|\ph|^2+|u_1|^2+|u_2|^2\bigr)\,.
\end{align} 
Now, we make the \pier{choice $\,M:=1/\bigl(4\chi^2\bigr)$.} Then the inequality
\eqref{cont:dep:weak} follows from an application of Gronwall's lemma. As a consequence, 
if $u_1^1=u_2^1$ and $u_1^2=u_2^2$, then  $\mu_1=\mu_2$, $\ph_1=\ph_2$ and $\s_1=\s_2$. But then, by \eqref{diff2}, also $\xi_1=\xi_2$. That is, the solution is  unique.   
The assertion is thus
completely proved.
\end{proof}

\section{\pier{Regularity properties}}
\label{REGUL}
\setcounter{equation}{0}
\betti{The next theorem provides a regularity result in \pier{the case of a general potential $F$ satisfying \ref{F:weak}, but} under more regular initial data and sources $(u_1,u_2)$  with respect to \pier{the well-posedness result in Theorem~\ref{THM:WEAK}.}}
\begin{theorem}[Regularity]
\label{THM:REGU}
\betti{Assume that \ref{const:weak}--\ref{h:weak} hold, and \pier{let} the initial data  satisfy \eqref{weak:initialdata} as well as the \pier{additional} assumptions
\begin{align}
	\label{strong:initialdata}
	\mu_0\in W,\quad \mu_0'\in V, \quad \s_0 \in L^\infty(\Omega), \quad (\partial F_1)^\circ (\ph_0)\in L^\infty(\Omega)\,.
\end{align}
Moreover, suppose that 
\begin{align}
	\label{u:strong}
		(u_1, u_2) \in L^2(0,T;V) \times \pier{L^\infty(0,T;H)}.
\end{align}
Then} \pier{the solution $(\m,\ph,\xi,\s)$ to \Statesys\ in the sense of Definition \ref{DEF:WEAK}
enjoys the further regularity properties}
\begin{equation}
\label{betti1}
\pier{\sigma\in L^\infty (Q), \quad \mu\in H^1(0,T; V)\cap L^\infty(0,T;W), \quad  \xi \in L^\infty(Q)\,.}
\end{equation}
\end{theorem}
\begin{proof}
First, observe that condition~\eqref{u:strong} implies~\eqref{u:uniq} by Sobolev embeddings. Therefore, Theorem~\ref{THM:WEAK} already ensures the existence of a unique solution. Consequently, as in the proof of the existence result, we proceed in a formal manner, employing the Yosida approximation of~$\partial F_{1}$ in our estimates, without resorting to finite-dimensional approximation techniques.

As in the previous proof, and for ease of notation, we continue to denote by $(\mu, \varphi, \xi, \sigma)$ the solution of the \revis{$\varepsilon$-approximating} system, where $\xi = F'_{1,\varepsilon}(\varphi)$, throughout the computations below. We will revert to the notation $(\mu_\varepsilon, \varphi_\varepsilon, \xi_\varepsilon, \sigma_\varepsilon)$ at the end of each estimate.

\noindent
{\sc First estimate.} From the boundedness properties in \eqref{smokeyuno} and \eqref{smokeydue}, along with~\eqref{u:strong}, we infer that the right-hand side of~\eqref{var:3} is uniformly bounded in $L^\infty(0,T;H)$. Since the initial datum $\sigma_0$ belongs to $L^\infty(\Omega)$ (cf.~\eqref{strong:initialdata}), maximal parabolic regularity (see, e.g., \cite[Chapter~III, Theorem~7.1, p.~181]{Ladyzenskaja_69}) yields
\begin{equation}\label{reg:sig}
  \|\sigma_\varepsilon\|_{L^\infty(Q)} \leq C .
\end{equation}

\noindent
{\sc Second estimate.}  We test \eqref{var:1} by $-\Delta \mu_t$ and integrate over $(0,t)$. Integrating by parts in space, and using the assumptions~\ref{P:weak}, \ref{h:weak}, \eqref{strong:initialdata}, \eqref{u:strong} along with the estimates 
\eqref{smokeyuno}, \eqref{smokeydue}, \juerg{we obtain from  the H\"older and Young inequalities that} 
\begin{align}
&\frac{\alpha}{2} \|\nabla\dt \mu(t)\|^2+\frac{1}{2} \|\Delta \mu (t)\|^2\non\\
&{}=\frac \a 2 \|\mu_0'\|_V^2+\frac 12 \|\Delta \mu_0\|^2 - \int_{Q_t}\nabla \dt\ph\cdot \nabla\dt\mu \non\\
&\quad{}+ \int_{Q_t}\bigl[P'(\ph)(\sigma+\chi(1-\ph)\jupi{{}-\mu{}}) - u_1\h'(\ph)\bigr]
\nabla \ph \cdot \nabla\dt 
\mu\non\\
&\quad{} +\int_{Q_t}P(\ph)\nabla (\sigma-\chi\ph\jupi{{}-\mu{}})\cdot \nabla\dt\mu
-\int_{Q_t}\h (\ph)\nabla u_1 \cdot \nabla\dt\mu\non\\
\separa
&{}\leq C + C \int_{Q_t} |\nabla \dt\ph|^2 
\non\\
&\quad{}
+ C \int_0^t \norma{\bigl( |\sigma| + 1 + |\ph| + | \mu| + |u_1| \bigr)(s) }^2_{3} 
\norma{\nabla \ph(s) }^2_{6}\, 
ds
\non\\
&\quad{}
+C \int_{Q_t} |\nabla (\sigma-\chi\ph \jupi{\,-\,} \mu)|^2 +C \int_{Q_t} |\nabla u_1|^2  
+ \int_0^t \|\nabla\dt\mu\|^2\,. \non
\end{align}
Note that all terms on the \rhs\ except the last are already bounded due to \eqref{smokeyuno}, \eqref{smokeydue}, \eqref{u:strong} and Sobolev embeddings. Then, 
applying now a standard Gronwall lemma, we derive the regularity estimate 
\begin{equation}\label{reg:mu}
\|\mu_\eps \|_{H^1(0,T; V)\cap L^\infty(0,T;W)}\leq C\,. 
\end{equation}

\noindent
{\sc Third estimate.} 
In view of \eqref{reg:mu}, it turns out that $\mu$ is bounded in $\L\infty W$,
hence in $\LQ\infty$. The same can be concluded for $F_2^{\,\prime}(\ph)$, due to \eqref{pier2-1} and \ref{F:weak}. 
Then, let us rewrite the equation~\eqref{pier6} as 
\begin{align}
  &\tau \dt \ph -\Delta \ph + F^{\,\prime}_{1,\eps}(\ph)  = h
  \quad \hbox{a.e. in }\,Q ,\non \\
  &\quad \hbox{with }\, h = \mu + \chi \,\sigma - F^{\,\prime}_{2}(\ph) \, \hbox{ bounded in } \,
  L^\infty (Q). 
  \label{seconda-mp}
\end{align}
To prove the third property in \eqref{betti1}, it is enough to derive {a} uniform $\LQ\infty$-bound for~$F^{\,\prime}_{1,\eps}(\ph)$.
Let us outline the argument by proceeding formally {and}
pointing out that {just} a truncation of {the} test functions would be needed for a rigorous proof.
We take any $p>2$ and test \eqref{seconda-mp} by $|F^{\,\prime}_{1,\eps}(\ph)|^{p-2}F^{\,\prime}_{1,\eps}(\ph)$, a function of $\ph$ which is increasing and attains the value 
$0$ at $0$ (cf.~\eqref{pier4}). Then, we integrate from $0$  to $t\in (0,T]$, obtaining
\begin{align}
  &{\tau}\iO \Bigl( \int_0^{\ph(t)} |F^{\,\prime}_{1,\eps}(s)|^{p-2}F^{\,\prime}_{1,\eps}(s)ds  \Bigr)
  \non 
  \\
  &{}
 + (p-1) \revis{\int_{Q_t}} |F^{\,\prime}_{1,\eps}(\ph)|^{p-2}F_{1,\eps}^{\,\prime\prime}(\ph) |\nabla\ph|^2
  +\revis{\int_{Q_t}} |F^{\,\prime}_{1,\eps}(\ph)|^p
  \non 
  \\
  &{}= {\tau} \iO \Bigl( \int_0^{\ph_0} |F^{\,\prime}_{1,\eps}(s)|^{p-2}F^{\,\prime}_{1,\eps}(s)ds \Bigr) + 
  \revis{\int_{Q_t}} h |F^{\,\prime}_{1,\eps}(\ph)|^{p-2}F^{\,\prime}_{1,\eps}(\ph) \,.
  \label{pier9}
\end{align}
Note that the first term {on the \lhs\ is nonnegative since $F^{\,\prime}_{1,\eps}$ is monotone increasing with $F^{\,\prime}_{1,\eps}(0)=0$;
moreover, the second term on the \lhs\ is nonnegative} since the derivative $F^{\,\prime\prime}_{1,\eps}$ is 
nonnegative everywhere.  About the \rhs\ we may recall \eqref{pier5} and observe that 
$$ {\tau} \iO \Bigl( \int_0^{\ph_0} |F^{\,\prime}_{1,\eps}(s)|^{p-2}F^{\,\prime}_{1,\eps}(s)ds \Bigr) \leq {\tau} \norma{ F^{\,\prime}_{1,\eps}(\ph_0)}_\infty^{p-1} 
\norma{ \ph_0}_\infty |\Omega|\,\leq C \norma{ (\partial F_1)^\circ (\ph_0)}_\infty^{p-1} . $$ 
On the other hand, thanks to {$p'=p/(p-1)$ and} the Young inequality, we have that
\begin{align}
  &\revis{\int_{Q_t}} h |F^{\,\prime}_{1,\eps}(\ph)|^{p-2}F^{\,\prime}_{1,\eps}(\ph)
  \leq \norma h_{L^p (Q_t)} \, \norma{\,|F^{\,\prime}_{1,\eps}(\ph)|^{p-1}}_{L^{p'} (Q_t)}
  \non
  \\
  &\quad{}
  = \norma h_{L^p (Q_t)} \, \norma{F^{\,\prime}_{1,\eps}(\ph)}_{L^p (Q_t)}^{p/p'} 
  \leq \frac 1p \, \norma h_{L^p (Q_t)}^p
  + \frac 1{p'} \, \norma{F^{\,\prime}_{1,\eps}(\ph)}_{L^p (Q_t)}^p \,.
  \non
\end{align}
By rearranging {from} \eqref{pier9}, and taking $t=T$, we infer that 
\begin{align}
  & \norma{F^{\,\prime}_{1,\eps}(\ph)}_{L^p (Q)} \leq  \Bigl( p \,C \norma{ (\partial F_1)^\circ (\ph_0)}_\infty^{p-1}+  \norma h_{L^p (Q)}^p \Bigr)^{1/p}  
  \non
  \\
  &\quad{}
  \leq  \bigl( p  \,C \norma{ (\partial F_1)^\circ (\ph_0)}_\infty^{p-1} \bigr)^{1/p}  +  \norma h_{L^p (Q)}.
\non
\end{align}
Then, {passing to the limit as $p\to +\infty$ in the above chain of inequalities}, we conclude~that 
$$\norma{F^{\,\prime}_{1,\eps}(\ph)}_{{L^\infty} (Q)} 
\leq  \norma{ \revis{(\partial F_1)^\circ (\ph_0)}}_\infty
  +  \norma h_{{L^\infty} (Q)}. $$
Hence, as $\xi_\eps = F^{\,\prime}_{1,\eps}(\ph_\eps)$ (now using the notation with dependence on $\eps$), \juerg{we finally arrive}~at 
\begin{equation}\label{reg:xi}
	\|\xi_\eps\|_{L^{\infty}(Q)} \leq C .
\end{equation}
We now collect the estimates \eqref{reg:sig}, \eqref{reg:mu}, and \eqref{reg:xi}, and recall the limiting procedure as $\varepsilon \searrow 0$, already carried out in the proof of Theorem~\ref{THM:WEAK}. Since the uniform estimates are preserved by lower semicontinuity, the proof is thus concluded.
\end{proof}

\Brem
\label{separation}
Let us point out that the regularity $\xi \in L^\infty(Q)$ established in~\eqref{betti1} implies, in particular, the so-called \emph{separation property} in the case of the logarithmic potential~\eqref{logpot}. Indeed, since
\[
  \xi = \ln \frac{1+\ph}{1-\ph} \in L^\infty(Q),
\]
there exist two real numbers $r_*$ and $r^*$, depending on the structure of the system, such that
\begin{align}
  -1 < r_* \leq \ph(x,t) \leq r^* < 1 
  \quad \text{for every } (x,t) \in \overline{Q}.
  \label{stima-sep}
\end{align}
Hence, the boundedness of $\xi$ prevents the phase variable $\ph$ from approaching the singular values $\pm 1$, ensuring that $\ph$ remains within physically meaningful bounds throughout the evolution.
\Erem

\section{The case \juerg{when $P$ is} constant}
\label{Asy}
\setcounter{equation}{0}

\def\muza{\mu_{0, \alpha}}
\def\nuza{\mu'_{0, \alpha}}

\def\mua{\mu_{\alpha}}

\def\sigza{\sigma_{0, \alpha}}

In this section, we restrict ourselves to the case where the proliferation function $P$ reduces to a positive constant and investigate the asymptotic behavior of the problem as $\alpha \searrow 0$. 
Accordingly, we strengthen assumption~\ref{P:weak} by imposing
\begin{equation}
\label{Pconst}
P \text{ is a positive constant.}
\end{equation}
Besides, we allow the initial data for $\mu$, $\partial_t \mu$, and $\sigma$ to depend on $\alpha$, while we keep $\ph_0$, the initial value of $\ph$, fixed. This choice is made for simplicity, in view of the restrictions on $\ph_0$ stated in~\eqref{weak:initialdata}.

Thus, for $0< \alpha \leq 1$, we consider families of data $ \{\mu_{0,\alpha} ,  \mu'_{0,\alpha}, \sigza \}$ such that 
\begin{align}
& \mu_{0,\alpha}\, \hbox{ is uniformly bounded in }\,  V ,\label{hyp-p1}\\[1mm]
& \mu'_{0,\alpha} \, \hbox{ is uniformly bounded in }\,  H , \label{hyp-p2}\\[1mm]
&  \sigza  \, \hbox{ is uniformly bounded in }\,  V \non \\
& \qquad \hbox{and strongly converges to $\sigma_0$ in } \, H \, \hbox{ as }
\alpha\searrow 0. \label{hyp-p3}
\end{align} 
Of course, it follows from \eqref{hyp-p3} that $\sigma_0 \in V $ (cf.~\eqref{weak:initialdata}) and $\sigza \to \sigma_0$ weakly in $V$. 
We can state the following convergence result. 

\Bthm
\label{Convergence}
Assume that \ref{const:weak}--\ref{h:weak}, \eqref{Pconst} hold, and \pier{let} the initial data satisfy \eqref{weak:initialdata} and \eqref{hyp-p1}--\eqref{hyp-p3}. Moreover, let 
\begin{align}
	\label{u:asym}
		(u_1, u_2) \in L^2(Q) \times L^2(Q) \quad \hbox{with} \quad
		u_1 \in \H1H \cap \L2{\Lx3},
\end{align}
so that also \eqref{u:uniq} is satisfied. 
For all $\alpha \in (0,1]$, let the quadruple $(\mu_\alpha , \ph_\alpha , \xi_\alpha, \sigma_\alpha)$, with 
\begin{align}
  & \mu_\alpha \in \W{2,\infty}{V^*}\cap \W{1,\infty}H \cap \L\infty V,
  \label{regmual}
  \\
  & \ph_\alpha \in W^{1,\infty}(0,T;H)  \cap H^1(0,T;V) \cap L^\infty(0,T;W)\cap C^0(\overline Q),
  \label{regphial}
  \\
  & \xi_\alpha \in \L\infty H,
  \label{regxial}
  	\\
  &\s_\alpha \in \H1 H \cap \L\infty V \cap \L2 {W}, \label{regsal}
\end{align}
be the solution \juerg{to} the initial value problem
\begin{align}
	 \label{primaal}& \<\alpha\dtt\mu_\alpha , v > + \iO\dt\ph_\alpha \, v 
	+ \iO \nabla \mu_\alpha \cdot \nabla v
	\nonumber\\
	&{}= \iO P\, (\sigma_\alpha+\bettib{\chi}(1-\ph_\alpha)-\mu_\alpha)v
	-\iO \h(\ph_\alpha)u_1 v
	 \nonumber\\
	&\qquad{} \hbox{for every $v \in V $ and a.e. in $(0,T)$,}\\[2mm]
	\label{secondaal}&\pier{\tau}\dt\vp_\alpha-\Delta\vp_\alpha+\xi_\alpha+F_2^{\,\prime}(\vp_\alpha)=\mu_\alpha+\chi\,\sigma_\alpha, \quad \hbox{$\xi_\alpha \in \partial F_1(\ph_\alpha)$, \, a.e. in $\,Q$,}
 	\\[2mm]
 	\label{terzaal}
 	& \dt\sigma_\alpha -\Delta\sigma_\alpha
 	=-\chi\Delta\ph_\alpha 
 	-  P\, (\sigma_\alpha+\bettib{\chi}(1-\ph_\alpha)-\mu_\alpha)
 	+ u_2	 \, \,\,\mbox{ a.e. in \,$Q$,}
	\\[2mm]
	\label{cauchyal}
	&\mua(0)=\muza \, ,\quad  (\dt \mua) (0)=\nuza \, ,\quad  \ph_\alpha(0)= \ph_0, \quad \s(0)=\sigza \quad \hbox{a.e. in $\Omega$}.
\end{align}
Then there exists a quadruple $(\mu,\ph,\xi, \s)$ such that,
for some subsequence $\alpha_k$ tending to $0$, there holds  
\begin{align}
\label{conv7}
&\mu_{\alpha_k} \to\mu \quad \mbox{weakly star in }\,\L\infty V,\\
\label{conv8}
&\alpha_k\mu_{\alpha_k} \to 0 \quad \mbox{weakly star in }\,\W{2,\infty}{V^*} \,  
\mbox{ and strongly in }\,\W{1,\infty}H,\\
\label{conv9}
&\ph_{\alpha_k}\to\ph \quad \mbox{weakly star in } \,W^{1,\infty}(0,T;H)\cap H^1(0,T;V)\cap L^\infty(0,T;W)
\non
\\
&\qquad 
\mbox{ and strongly in }\,\C0V \cap C^0(\overline Q),\\
\label{conv10}
&\xi_{\alpha_k}\to\xi \quad \mbox{weakly star in } L^\infty(0,T;H),
\\
\label{conv11}
&\s_{\alpha_k}\to\s \quad \mbox{weakly star in } 
\H1 H \cap \L\infty V \cap \L2 {W}.
\end{align}
Moreover, $(\mu,\ph,\xi, \s)$ is a solution \juerg{to} the viscous Cahn--Hilliard system 
\begin{align}
	 \label{pier:1}& \iO\dt\ph \, v 
	+ \iO \nabla \mu \cdot \nabla v
	= \iO P\,(\sigma+\chi(1-\ph)-\mu)v
	-\iO \h(\ph)u_1 v
	\nonumber\\
	& \qquad \hbox{for every $v \in V $ and a.e. in $(0,T)$,}\\[2mm]
	\label{pier:2}&\pier{\tau}\dt\vp-\Delta\vp+\xi+F_2^{\,\prime}(\vp)=\mu+\chi\,\sigma, \quad \hbox{$\xi \in \partial F_1(\ph)$, \, a.e. in $\,Q$,}
 	\\
 	\label{pier:3}
 	& \dt\sigma -\Delta\sigma
 	=-\chi\Delta\ph 
 	-  P\, (\sigma+\chi(1-\ph)-\mu)
 	+ u_2	 \, \,\,\mbox{ a.e. in \,$Q$,}\\
	\label{pier:4}
	&\ph(0)=\ph_0, \quad
	\s(0)=\s_0, \quad \hbox{a.e. in $\Omega$}.
\end{align}
\Ethm

\Bdim
We go back to the proof of Theorem~\ref{THM:WEAK} and consider the First Estimate. We focus, in particular, on  the equality~\eqref{esti1}. Now, we change the treatment of the term coming from the \rhs\ of~\eqref{primaal}. Please, let us use the notation without any index when doing the computation. Then, by integrating by parts in time, in place of~\eqref{esti1pier} we obtain  
\begin{align}
&\int_{Q_t}\bigl[P\, (\s+\jupi{\chi (1-\ph)}-\mu)-\h(\ph)u_1\bigr]\dt\mu
\non\\
&={} - \frac P 2\iO |\mu(t)|^2  +  \frac P 2\iO |\mu_{0,\alpha}|^2 + 
\iO [P\, (\s+\jupi{\chi (1-\ph)})-\h(\ph)u_1](t)\,  \mu(t) \non\\
&\quad{}
- \iO \bigl[P\, (\s_{0,\alpha} +\jupi{\chi (1-\ph_0)}) -\h(\ph_0) u_1(0) \bigr] \mu_{0,\alpha}
\non\\
&\quad{}
- \int_{Q_t} [P\, (\dt\s - \chi \, \dt \ph)-\h'(\ph)\dt \ph \, u_1 - \h (\ph) \dt u_1 ]\mu , \non
\end{align}
whence,  from \ref{h:weak}, \eqref{Pconst} and H\"older's inequality, it follows that
\begin{align}
&\int_{Q_t}\bigl[P\, (\s+\jupi{\chi (1-\ph)}-\mu)-\h(\ph)u_1\bigr]\dt\mu \non \\
&{}\le{} - \frac P 2\iO |\mu(t)|^2  + C
+ \frac P 4\iO |\mu(t)|^2 +  C\iO \bigl(|\s(t) |^2 \jupi{{}+1{}} +|\ph(t) |^2 \bigr) \non\\
&\quad{}+ C  \norma{u_1}_{\L\infty H}^2
+ C \bigl( \norma{\s_{0,\alpha}} \jupi{{}+1{}}+ \norma {\ph_0}  +  \norma{u_1}_{\L\infty H} \bigr) \norma{\mu_{0,\alpha}}
\non\\
&\quad{}
+\frac1{16}\int_{Q_t} |\dt\s|^2 +  C\int_{Q_t} \bigl( |\dt\ph|^2+ |\mu|^2 \bigr) 
\non\\
&\quad{}+ C \int_0^t \norma{\dt\ph(s) }
\,\norma{u_1(s) }_3 \, \norma{\mu (s)}_6 ds 
+ C \norma{\dt u_1}_{\L2H}^2  +  C\int_{Q_t} |\mu|^2\, . \non
\end{align}
Therefore, in view of \eqref{hyp-p1}--\eqref{u:asym}, and using Sobolev's embeddings and Young's inequality,  we have that
\begin{align*}
&\int_{Q_t}\bigl[P\, (\s+\jupi{\chi (1-\ph)}-\mu)-\h(\ph)u_1\bigr]\dt\mu\\
&{}\le{} - \frac P 4\iO |\mu(t)|^2  + C + C \left(\|\s_{0\alpha}\|^2+\intQt\s\dt\s\right) +C\left(\|\ph_0\|^2+\intQt\ph\dt\ph\right)\\
&\quad{}
+\frac1{16}\int_{Q_t} |\dt\s|^2 +  C\int_{Q_t} \bigl( |\dt\ph|^2+ |\mu|^2 \bigr) + C \int_0^t
\,\norma{u_1(s) }^2_3 \, \norma{\mu (s)}^2_V ds ,
\end{align*}
and, consequently,
\begin{align}
\label{pier21}
&\int_{Q_t}\bigl[P\, (\s+\jupi{\chi (1-\ph)}-\mu)-\h(\ph)u_1\bigr]\dt\mu \non \\
&{}\le - \frac P 4\iO |\mu(t)|^2  + C + \revis{\frac2{16}}\int_{Q_t} |\dt\s|^2
\non\\
&\quad{} +  C\int_{Q_t} \bigl(|\s|^2 +|\ph|^2 +  |\dt\ph|^2+ |\mu|^2 \bigr) + C \int_0^t
\,\norma{u_1(s) }^2_3 \, \norma{\mu (s)}^2_V ds 
\,.
\end{align}
By virtue of \eqref{pier21}, the inequality \eqref{esti3} now becomes
\begin{align}
\label{pier22}
&\frac\a 2\|\dt\mu(t)\|^2 +  \frac P 4\iO |\mu(t)|^2+\frac 12\iO |\nabla \mu(t)|^2 +\frac\tau 2\|\dt\ph(t)\|^2
+\intQt |\nabla\dt\ph|^2
\non\\
&\le \,C+C\intQt \bigl(|\mu|^2 + |\ph|^2   + |\s|^2 + |\dt\ph|^2
\bigr) 
\non\\
&\quad{}+{\frac 3{16}}\intQt |\dt\s|^2+ C \int_0^t
\,\norma{u_1(s) }^2_3 \, \norma{\mu (s)}^2_V ds \,,
\end{align}
so that in this context the subsequent inequality \eqref{esti6pier} reads 
\begin{align}
\label{pier23}
&\frac\a 2\|\dt\mu(t)\|^2 +\min\left\{\frac P 4, \frac 1 2\right\} \| \mu(t)\|_V^2 +\frac\tau 2\|\dt\ph(t)\|^2 +\frac{\chi^2}2\|\ph(t)\|_V^2 \non \\
&\quad +  4 \tau \chi^2 \intQt |\dt\ph|^2 + \frac12 \intQt |\nabla\dt\ph|^2 
+ \frac14 \intQt|\dt\s|^2 + \frac 14 \|\s(t)\|_V^2
\non\\
&\le \,C + C\intQt \bigl(|\mu|^2 + |\ph|^2   + |\s|^2 + |\dt\ph|^2
\bigr) + C \int_0^t
\,\norma{u_1(s) }^2_3 \, \norma{\mu (s)}^2_V ds  \,,
\end{align}
where all \juerg{of} the constants $C $ in the above sequence of estimates are independent of $\alpha$. Note that the function $s \mapsto \norma{u_1(s) }^2_3$ belongs to $L^1(0,T) $ (see~\eqref{u:asym}), whence the application of the Gronwall lemma to \eqref{pier23} leads us to 
\begin{align}
\label{pier24}
&\alpha^{1/2}\|\mu\|_{W^{1,\infty}(0,T;H)}\,+\,\|\mu\|_{L^\infty(0,T;V)} \,+\,
\|\ph\|_{W^{1,\infty}(0,T;H)\cap H^1(0,T;V)}\non\\
&\quad + \|\s\|_{H^1(0,T;H)\cap L^\infty(0,T;V)}\,\le\,C \, ,
\end{align}
which replaces \eqref{smokeyuno}. 
Next, a closer inspection of the proof of Theorem~\ref{THM:WEAK}, Second Estimate,  reveals that the estimates~\eqref{smokeydue}--\eqref{smokey3bis} still hold,
with constants independent of $\alpha$, and that \eqref{smokeyfour} can be replaced by 
\begin{align}
\label{pier25}
&\alpha \|\mu\|_{W^{2,\infty}(0,T;V^*)}\,\le\,C\,,
\end{align}
on account of the fact that (cf.~\eqref{u:asym}) the \rhs\ of \eqref{primaal} is now bounded in~$ \L\infty H$. In conclusion, the boundedness properties \eqref{smokeydue}--\eqref{smokey3bis}, \eqref{pier24}, \eqref{pier25} are also valid for the solution $(\mu_\alpha , \ph_\alpha , \xi_\alpha, \sigma_\alpha)$ of \eqref{primaal}--\eqref{cauchyal}. This solution is unique due to Theorem~\ref{THM:WEAK}, since \eqref{u:asym} implies \eqref{u:uniq}. 

Then, by a 
standard weak star compactness argument, we deduce the existence of a subsequence $\alpha_k\searrow 0$ 
and a quadruple~$(\mu,\ph,\xi, \s)$ such that \eqref{conv7}--\eqref{conv11} hold. 
At this point, we can perform the limit procedure as in the passage to the limit as $\eps \searrow 0$ in the proof of Theorem~\ref{THM:WEAK}. We just point out that\juerg{, in order} to obtain the inclusion in \eqref{pier:2}, we should instead use~\cite[{Lemma~2.3,} p.~38]{Barbu}, since the subdifferential structure explicitly appears in~\eqref{secondaal} as well. This concludes the proof.
\Edim

\Brem
\label{lim-prob}
Note that Theorem~\ref{Convergence} implicitly guarantees the existence of
solutions to system~\eqref{pier:1}--\eqref{pier:4} for all potentials $F$
satisfying~\ref{F:weak}, and therefore for every convex and lower
semicontinuous function $F_1:\mathbb{R}\to [0,+\infty]$ with $F_1(0)=0$.
The obtained solution is in fact a strong solution. Indeed, using
\eqref{conv7}--\eqref{conv11}, treating \eqref{pier:1} appropriately, and
invoking elliptic regularity theory, we see that \eqref{pier:1} can be rewritten as
\begin{align}
   &\partial_t \varphi - \Delta \mu
      = P\,(\sigma+\chi(1-\varphi)-\mu)
        - \revis{\h}(\varphi)u_1
        \qquad \text{a.e.\ in } Q,
   \label{pier14}
\end{align}
supplemented with the boundary condition
\begin{equation}
   \partial_{\mathbf{n}} \mu = 0
   \qquad \text{a.e.\ on } \Sigma ,
   \label{pier16}
\end{equation}
which yields the additional regularity $\mu \in L^\infty(0,T;W)$.
Moreover, we remark that if $\sigma_0 \in V \cap L^\infty(\Omega)$ and
$u_2 \in L^\infty(0,T;H)$, then, by arguing as in the First Estimate in the
proof of Theorem~\ref{THM:REGU}, we obtain the regularity
$\sigma \in L^\infty(Q)$.

At this stage, if $(\partial F_1)^\circ(\varphi_0)\in L^\infty(\Omega)$,
we may also repeat the Third Estimate from the proof of
Theorem~\ref{THM:REGU}, since the right-hand side $h$ in
\eqref{seconda-mp} belongs to $L^\infty(Q)$. Hence, we deduce that
$\xi \in L^\infty(Q)$, which in particular yields a separation property
(cf.\ Remark~\ref{separation}) in the case of the logarithmic
potential~\eqref{logpot}.
\Erem

\def\mubar{\overline{\mu}}
\def\phbar{\overline{\ph}}
\def\xibar{\overline{\xi}}
\def\sbar{\overline{\s}}
\def\rhobar{\overline{\rho}}
\def\szerobar{\overline{\s}_0}

\pcol{We now derive an error estimate for the differences $\varphi_\alpha - \varphi$ and $\sigma_\alpha - \sigma$, measured in suitable norms and quantified in terms of the parameter $\alpha$.}

\Bthm
\label{Errest}
\pcol{Under the same assumptions as in Theorem~\ref{Convergence}, we let $(\mu_\alpha , \ph_\alpha , \xi_\alpha, \betti{\s}_\alpha)$ denote the solution to \eqref{primaal}--\eqref{cauchyal}, for $\alpha \in (0,
1]$, and $(\mu,\ph,\xi, \s)$ be the solution to \eqref{pier:1}--\eqref{pier:4} found by the 
asymptotic limit in \eqref{conv7}--\eqref{conv11}. Moreover, in addition to~\eqref{hyp-p3} we assume as well that 
\begin{align}
&  \norma{\sigza- \sigma_0} \leq C_\sigma \, \alpha^{1/4} \quad \hbox{ for every $\alpha \in (0,1] $,}
\label{hyp-p4}
\end{align} 
for some constant $C_\sigma>0$.
Then there is a constant $K_2>0$, which depends on the structure of the system but is independent of $\alpha$, such that}
\begin{align}
&\pcol{\alpha^{1/2} \norma{\mu_a}_{\L\infty H}
+ \norma{1* (\mua-\mu)}_{\L\infty V}
+ \norma{\ph_\alpha-\ph}_{\L\infty H\cap\L2V}}
\non \\
&\pcol{{}+ \norma{\s_\alpha-\s}_{\L2 H} 
+ \norma{1*(\s_\alpha-\s) }_{\L\infty V}
\leq K_2 \, \alpha^{1/4} \quad \hbox{ for every $\, \alpha \in (0,1] $.}}
  \label{stimaK2}
\end{align}
\Ethm

\Bdim
We proceed partly as in the proof of Theorem~\ref{THM:WEAK}, specifically for the steps concerning Uniqueness and Continuous Dependence. Let us introduce the auxiliary elements  $\rho_\alpha := 
\s_\alpha -\chi \ph_\alpha$, $\, \rho := \s -\chi \ph$ and set
\begin{align}
&\mubar := \mu_\alpha-\mu , \quad \phbar :=\ph_\alpha-\ph ,\quad \xibar :=\xi_\alpha -\xi ,\non\\
&\sbar := \s_\alpha-\s , \quad \rhobar := \rho_\alpha-\rho \quad \hbox{and} \quad 
\szerobar := \s_{0,\alpha} -\s_0. \non
\end{align}
The plan is to subtract \eqref{pier:1}–\eqref{pier:3} from the corresponding relations \eqref{primaal}–\eqref{terzaal}. We begin by integrating in time the difference between \eqref{primaal} and \eqref{pier:1}. Using \eqref{cauchyal}, we then deduce that
\begin{align}
 \label{diff-p1}
& \iO \alpha \, \dt \mua  \, v 
	+ \iO \nabla (1*\mubar) \cdot \nabla v + \iO P (1*\mubar) v  = \iO  \alpha \, \nuza \, v - \iO \phbar\, v
\non\\
&\quad{} + \iO P ( 1* \rhobar)v
-\iO [1* ((\h(\ph_\alpha)-\h(\ph ))u_1)] v\nonumber\\
	& \qquad{} \hbox{for every $v \in V $ and a.e. in $(0,T)$.}
\end{align}
From \eqref{secondaal} and \eqref{pier:2} it follows that
\begin{align}
\label{diff-p2}
&\pier{\tau}\dt\phbar-\Delta\phbar+\xibar +F_2^{\,\prime}(\vp_\alpha)-F_2^{\,\prime}(\ph )=\mubar +\chi\,\rhobar + \chi^2 \, \phbar, 
\non\\
&\quad{} \hbox{$\xi_\alpha \in \partial F_1(\ph_\alpha)$, \ $\xi \in \partial F_1(\ph)$, \ a.e. in $\,Q$.}
\end{align}
We then perform another time integration on the difference of \eqref{terzaal} and \eqref{pier:3}, using \eqref{cauchyal} together with \eqref{pier:4}.
After adding $-\chi\,\bettib{\phbar}$ to both sides, it is straightforward to arrive~at
\begin{align}
\label{diff-p3}
 	& \rhobar -\Delta (1*\rhobar) + P (1* \rhobar)
 	= \szerobar +  P (1*\mubar ) - \chi\, \bettib{\phbar}
  	 \, \,\,\mbox{ a.e. in \,$Q$.}
\end{align}
We also note that (cf.~\eqref{convolution})
\begin{align} 
\label{diff-p4}
&(1*\mubar)(0)=0, \quad \bettib{\phbar}\/(0)=0,\quad (1*\rhobar) (0)=0,\quad
\mbox{ a.e. in $\,\Omega$}. 
\end{align}

Now we take $v = \mubar$ in \eqref{diff-p1} and test the equality in \eqref{diff-p2} 
by $\phbar$. Summing the two relations, observing a cancellation, and 
integrating with respect to time, we obtain the inequality below.  
Since the product $\xibar \, \phbar$ is nonnegative -- by the inclusions in 
\eqref{diff-p2} and the monotonicity of $\partial F_1$ -- we easily derive
\begin{align}
&\frac \alpha 2\| \mua (t)\|^2 + \frac{\bettib{K}}{ 2 }\norma{(1*\mubar) (t)}_V^2 + 
\frac \tau 2\| \phbar (t)\|^2 
+ \intQt |\nabla \phbar |^2 
\non\\
& \leq\,
\frac \alpha 2\|\muza \|^2 + \intQt \alpha \, \dt \mua \, \mu  +
\iO  \alpha \, \nuza \pcol{(1*\mubar)(t)}  
\non\\
&\quad{}+ \intQt P ( 1* \rhobar)\, \bettib{\mubar}  -\intQt [1* ((\h(\ph_\alpha)-\h(\ph ))u_1)]\, \bettib{\mubar} 
\non\\
&\quad{} - \intQt \bigl(F_2^{\,\prime} (\ph_\alpha) -  F_2^{\,\prime} (\ph) - \chi^2\, \phbar \bigr) \phbar + \intQt \chi\, \rhobar \, \phbar
\label{pier19}
\end{align}
for all $t\in (0,T]$, where $\bettib{K}:= \min  \{P,1\}$. Now, recalling the boundedness 
properties~\eqref{hyp-p1} and \eqref{hyp-p2}, the estimate  \eqref{pier24} for 
$\alpha^{1/2} \norma{\dt\mua}_{\L\infty H} $, as well as the regularity  $\mu \in \L\infty H$, \betti{we deduce from Young's inequality}
\begin{align}
&\frac \alpha 2\|\muza \|^2 + \intQt \alpha \, \dt \mua \, \mu  +
\iO  \alpha \, \nuza \pcol{(1*\mubar)(t)} 
\non\\
&{}\leq 
C \alpha  + {C} \alpha^{1/2}  + \frac {\bettib{K}} {16} \iO|(1*\mubar) (t)|^2 + C \alpha^2 \norma{\nuza}^2
\non\\
&{}\leq C \alpha^{1/2} + \frac {\bettib{K}} {16} \iO|(1*\mubar) (t)|^2 .
\label{pier26}
\end{align}
We integrate by parts in time the next two terms on the right side of \eqref{pier19}. With the aid  of the Lipschitz continuity of $\h$ (see~\eqref{differ1}), the estimates \eqref{pieruwe} and \eqref{uwe} with $p=6$, and H\"older and Young's inequalities we obtain 
\begin{align}
& \intQt P ( 1* \rhobar) \,\bettib{\mubar} -\intQt [1* ((\h(\ph_\alpha)-\h(\ph ))u_1)]\, \bettib{\mubar} 
\non\\
&= \iO P ( 1* \rhobar)(t) (1* \bettib{\mubar})(t)  - \iO [1* ((\h(\ph_\alpha)-\h(\ph ))u_1)](t)  (1* \bettib{\mubar})(t)
\non\\
&\quad{}- \intQt P \, \rhobar\, ( 1* \mubar) + \intQt  (\h(\ph_\alpha)-\h(\ph ))u_1 ( 1* \mubar)
\non\\
&\leq \frac{\bettib{K}} {16} \iO|(1*\mubar) (t)|^2 + \frac{4P^2}  {\bettib{K}} \iO|(1*\rhobar) (t)|^2
+ L  \norma{ (1*\pcol{\mubar})(t)}_6\int_0^t \norma{\phbar(s)}\norma{u_1(s) }_3\, ds\non\\
&\quad{}
+  \frac 1 4 \intQt |\rhobar|^2  +   {P^2}  \int_0^t \norma{ (1* \mubar)(s)}^2\,ds + 
L \int_0^t \norma{\phbar(s)}\norma{u_1(s) }_3\norma{ (1*\pcol{\mubar})(s)}_6 \, ds .
\non
\end{align}
Hence, thanks to the Sobolev embeddings and again to the H\"older and Young inequalities, we infer that 
\begin{align}
& \intQt P ( 1* \rhobar) \bettib{\mubar} -\intQt [1* ((\h(\ph_\alpha)-\h(\ph ))u_1)]\, \bettib{\mubar} 
\non\\
&\leq \frac{\bettib{K}} {8} \norma{(1*\mubar) (t)}_V^2 + \frac{4P^2}  {\bettib{K}} \iO|(1*\rhobar) (t)|^2
 \, ds +  \frac 1 4 \intQt |\rhobar|^2\non\\
&\quad{}
 +   C \int_0^t \norma{ (1* \mubar)(s)}_V^2\,ds +  C  \int_0^t \norma{u_1(s) }_3^2 \, \norma{\phbar(s)}^2 \, ds .
\label{pier27}
\end{align}
As for the last two terms of \eqref{pier19}, by virtue of the Lipschitz continuity of $F_2^{\,\prime}$ and the Young inequality, we have that 
\begin{align}
&- \intQt \bigl(F_2^{\,\prime} (\ph_\alpha) -  F_2^{\,\prime} (\ph) - \chi^2\, \phbar \bigr) \phbar + \intQt \chi\, \rhobar \, \phbar
\leq  \frac 1 4 \intQt |\rhobar|^2 + C \intQt |\phbar|^2 . 
\label{pier28}
\end{align}

At this point, collecting the estimates \eqref{pier26}--\eqref{pier28}, \betti{it follows
from \eqref{pier19}} that 
\begin{align}
&\frac{\alpha}{2}\| \mua(t)\|^2 
 + \frac{\bettib{K}}{4} \bigl\|(1*\mubar)(t)\bigr\|_V^2
 + \frac{\tau}{2}\| \phbar(t)\|^2 
 + \int_{Q_t} |\nabla \phbar|^2 
\non\\
&\le\,
C \alpha^{1/2} 
 + \frac{4P^2}{\bettib{K}} \int_\Omega |(1*\rhobar)(t)|^2 \, ds
 + \frac{1}{2} \int_{Q_t} |\rhobar|^2 
\non\\
&\quad{}
 + C \int_0^t \bigl\|(1*\mubar)(s)\bigr\|_V^2\, ds
 + C \int_0^t \bigl( 1 + \|u_1(s)\|_3^2 \bigr)\, \|\phbar(s)\|^2 \, ds
\label{pier29}
\end{align}
for all $t \in (0,T]$.  

Next, we test \eqref{diff-p3} by $M \rhobar$, where the constant $M>0$ will be specified below. 
Integrating with respect to time and applying Young's inequality, we obtain 
\begin{align}
& M \int_{Q_t} |\rhobar|^2 
  + \frac{M\bettib{K}}{2} \bigl\|(1*\rhobar)(t)\bigr\|_V^2
\le 
\frac{M}{2} \int_{Q_t} |\szerobar + P(1*\mubar) - \chi\,\bettib{\phbar}|^2
 + \frac{M}{2} \int_{Q_t} |\rhobar|^2,
\non
\end{align}
so that the assumption \eqref{hyp-p4} allows us to deduce that 
\begin{align}
& \frac{M}{2} \int_{Q_t} |\rhobar|^2 
  + \frac{M\bettib{K}}{2} \bigl\|(1*\rhobar)(t)\bigr\|_V^2
\le 
MC \alpha^{1/2}
 + MC \int_{Q_t} \bigl( |1*\mubar|^2 + |\bettib{\phbar}|^2 \bigr),
\label{pier30}
\end{align}
for all $t \in (0,T]$.  

Before adding \eqref{pier29} and \eqref{pier30}, we choose $M$ such that
\[
\frac{M}{2} > \frac{1}{2}
\quad\text{and}\quad
\frac{M\bettib{K}}{2} > \frac{4P^2}{\bettib{K}},
\qquad\text{i.e.},\qquad
M > \max\!\left\{ 1,\, \frac{8P^2}{\bettib{K}^2} \right\}.
\]
With this choice, we may sum \eqref{pier29} and \eqref{pier30} and then apply 
the Gronwall lemma, using the fact that 
\( s \mapsto \|u_1(s)\|_3^2 \) belongs to \( L^1(0,T) \) (cf.~\eqref{u:asym}). 
Therefore, we finally obtain the estimate 
\begin{align}
& \alpha^{1/2}\, \|\mu_\alpha\|_{L^\infty(0,T;H)}
  + \|1*\mubar\|_{L^\infty(0,T;V)}
  + \|\phbar\|_{L^\infty(0,T;H)\cap L^2(0,T;V)}
\non\\
&\qquad{}
  + \|\rhobar\|_{L^2(0,T;H)}
  + \|1*\rhobar\|_{L^\infty(0,T;V)}
  \le C\, \alpha^{1/4},
\label{pier31}
\end{align}
from which \eqref{stimaK2} follows immediately, recalling that 
\(
\rhobar = \rho_\alpha - \rho = \s_\alpha - \s - \chi(\ph_\alpha - \ph).
\)
\Edim

\Brem
\label{Unique}
We point out that applying a similar procedure to a pair of arbitrary solutions to problems \eqref{pier:1}--\eqref{pier:4} would allow us to prove the uniqueness of the quadruple $(\mu, \varphi, \xi, \sigma)$, with
\begin{align}
\mu &\in L^\infty(0,T;V), \qquad \xi \in L^\infty(0,T;H), \label{pier32}
\\
\varphi &\in W^{1,\infty}(0,T;H) \cap H^1(0,T;V) \cap L^\infty(0,T;W) \cap C^0(\overline Q), \label{pier33}
\\
\sigma &\in H^1(0,T;H) \cap C^0([0,T];V) \cap L^2(0,T;W), \label{pier34}
\end{align}
solving \eqref{pier:1}--\eqref{pier:4}. Indeed, obtaining an estimate analogous to \eqref{pier31} would directly yield uniqueness for $1*\mu$, $\varphi$, and $\rho$ (and consequently for $\sigma$), while the uniqueness of $\xi$ then follows from a comparison in \eqref{pier:2}. Moreover, the uniqueness property implies that the convergences \eqref{conv7}--\eqref{conv11} stated in Theorem~\ref{Convergence} hold not only along a subsequence $\alpha_k \to 0$, but for the entire family as $\alpha \to 0$.
\Erem


\section*{Acknowledgments}
\betti{PC and ER gratefully mention their affiliation
to the GNAMPA (Gruppo Nazionale per l'Analisi Matematica, 
la Probabilit\`a e le loro Applicazioni) of INdAM (Isti\-tuto 
Nazionale di Alta Matematica) as well as their collaboration,
as Research Associate, to the IMATI -- C.N.R. Pavia, Italy. 
PC and ER \pier{also acknowledge support from} the Next Generation EU Project 
No.~P2022Z7ZAJ (A unitary mathematical framework for modelling muscular dystrophies).}


\end{document}